\documentclass[12pt]{article} 
\usepackage{latexsym, 
amscd, 
amsthm, 
graphicx, epsfig, amssymb}

\newtheorem{thm}{Theorem}[section]
\newtheorem{defn}[thm]{Definition}
\newtheorem{prop}[thm]{Proposition}

\newtheorem{rema}[thm]{Remark}

\newtheorem{expl}[thm]{Example}

\newcommand{\halmos}{\rule{1ex}{1.4ex}}

\newcommand{\nn}{\nonumber \\}

 \newcommand{\res}{\mbox{\rm Res}}

\renewcommand{\hom}{\mbox{\rm Hom}}

 \newcommand{\pf}{{\it Proof.}\hspace{2ex}}
 
 \newcommand{\epfv}{\hspace*{\fill}\mbox{$\halmos$}\vspace{1em}}

\newcommand{\lbar}{\bigg\vert}

\newcommand{\C}{\mathbb{C}}
\newcommand{\Z}{\mathbb{Z}}

\newcommand{\N}{\mathbb{N}}

\newcommand{\one}{\mathbf{1}}

\title{ {\bf A cohomology theory of grading-restricted
vertex algebras} }
\date{}
\author{Yi-Zhi Huang}
\begin{document}

\bibliographystyle{alpha}
\maketitle

\begin{abstract}
We introduce a cohomology theory of grading-restricted
vertex algebras. To construct the {\it correct} cohomologies, 
we consider linear maps from tensor 
powers of a grading-restricted
vertex algebra to ``rational functions valued in the algebraic 
completion of a module for the algebra,'' instead of linear maps 
from tensor powers of the algebra to a module for the algebra. 
One subtle complication arising from such functions is that we have to  carefully
address the issue of convergence when we compose
these linear maps with vertex operators. In particular, for each $n\in \N$,
we have an inverse system $\{H^{n}_{m}(V, W)\}_{m\in \Z_{+}}$ of 
$n$-th cohomologies and an additional 
$n$-th cohomology $H_{\infty}^{n}(V, W)$
of a grading-restricted
vertex algebra $V$ with coefficients in a $V$-module $W$ 
such that $H_{\infty}^{n}(V, W)$ is isomorphic to the inverse limit
of the inverse system $\{H^{n}_{m}(V, W)\}_{m\in \Z_{+}}$. 
In the case of $n=2$, there is an additional 
second cohomology denoted by $H^{2}_{\frac{1}{2}}(V, W)$ 
which will be shown in
a sequel to the present paper 
to correspond to what we call 
square-zero extensions of $V$ and to first order 
deformations of $V$ when $W=V$. 
\end{abstract}

\renewcommand{\theequation}{\thesection.\arabic{equation}}
\renewcommand{\thethm}{\thesection.\arabic{thm}}
\setcounter{equation}{0}
\setcounter{thm}{0}

\section{Introduction}

Vertex (operator) algebras arose naturally in both mathematics and physics
(see \cite{BPZ}, \cite{B} and \cite{FLM}) and are analogous to both Lie algebras
and commutative associative algebras. In the studies of various algebraic structures,
including in particular Lie algebras, associative algebras, commutative associative
algebras, and their representations, the corresponding cohomology theories,
such as Chevalley-Eilenberg cohomology of Lie algebras \cite{CE}, Hochschild 
cohomology of associative algebras \cite{Ho} and Harrison  or Andr\'{e}-Quillen 
cohomology of commutative associative algebras \cite{Ha} \cite{A} \cite{Q},
play important roles. See, for example, the book \cite{W} for an excellent 
introduction to these theories (except for the Harrison cohomology). 
These cohomologies describe naturally
certain extensions of these algebras, extensions of their modules and 
also deformations of these algebras. 
Moreover, the powerful tool of homological algebra developed in the last sixty years
has been used to obtain many old and new results on these algebras and their modules. 
Though much progress, especially in the case of simple vertex operator
algebras satisfying certain finiteness and reductivity conditions, 
has been made in the theory of vertex (operator) algebras,
a {\it correct} cohomology theory of vertex (operator) algebras is urgently needed 
in order to have a better 
understanding of the structures of vertex (operator) algebras and their modules
and to use the powerful tool of homological algebra. 

In \cite{KV}, applying the general theory for algebras 
over operads developed by Ginzburg and Kapranov \cite{GK}, 
Kimura and Voronov introduced a cohomology theory of
algebras over the operad of the moduli
space of configurations of disjoint and ordered biholomorphic embeddings of the 
unit disk into the Riemann sphere. Motivated by the operadic and geometric formulation of 
vertex operator algebras by the author \cite{H1} \cite{H2} and by Lepowsky and the 
author \cite{HL1} \cite{HL2}, Kimura and Voronov proposed in \cite{KV} that 
their cohomology theory of algebras over the moduli space above
also gives the cohomology theory of vertex operator algebras. Unfortunately,
their proposal was based on the assumption that vertex operator algebras are
in particular algebras over the operad of the moduli space mentioned
above, while this assumption 
holds only for vertex operator 
algebras obtained from commutative associative algebras. In fact, given a vertex operator
algebra 
$$\biggl(V=\coprod_{n\in \Z}V_{(n)}, Y, \one, \omega\biggr)$$
that is not obtained from 
a commutative associative algebra, for $u, v\in V$ and $a, b, z\in \C^{\times}$, 
$$Y(a^{L(0)}u, z)b^{L(0)}v$$
is in general not an element of $V$, even when 
$a$, $b$ and $z$ are chosen such that they give 
the configuration of three disjoint and ordered
biholomorphic embeddings of the 
unit disk into the Riemann sphere. 
Instead, it is an element of the 
algebraic completion 
$$\overline{V}=\prod_{n\in \Z}V_{(n)}$$
of $V$. 
So for $z\in \C^{\times}$, the map 
$$Y(a^{L(0)}\cdot, z)b^{L(0)}\cdot: V\otimes V\to \overline{V}$$
in general does not belong to the 
endomorphism operad of the vector space $V$. Therefore, the vertex operator algebra 
$V$ in general does not give an algebra over the operad of the moduli space mentioned above. 
In particular, the cohomology theory introduced by Kimura and Voronov 
cannot be used to give a cohomology theory of vertex operator algebras 
because of this subtle but crucial feature of the geometric and operadic 
formulation of vertex operator algebra. Moreover, vertex operator algebras satisfy
an additional meromorphicity condition which one must take into account in
any cohomology theory of vertex operator algebras. 

In Section 11 of \cite{B2}, Borcherds also proposed a cohomology theory for general
vertex algebras by using his categorical formulation of vertex algebra and an
analogy with the Hochschild homology of associative algebras. 
However, the subtle details of this cohomology theory
were not carried out and the basic properties that
a cohomology theory must have were not discussed. 
More importantly, we are interested only in what
we call grading-restricted vertex algebras
and for these vertex algebras, a cohomology theory for general vertex algebras
cannot be the {\it correct} one. Here is the reason:
The notion of vertex algebra is too general to give
properties strong enough for a good representation theory.
The class of vertex algebras for which many substantial results in representation theory
have been obtained is
that of vertex operator algebras. 
Since we want to allow the deformations of 
the representation structures of the Virasoro 
algebra, especially the deformations of the central charges,
we are interested in the slightly 
more general class of  grading-restricted vertex algebras, for which 
conformal elements are not specified but $\Z$-gradings are still given and 
the grading-restriction condition is satisfied. 
For grading-restricted vertex algebras,
a cohomology 
theory for general vertex algebras cannot be the {\it correct} one,
because, for example, starting from a grading-restricted vertex algebra,
the deformations corresponding to such a general cohomology theory
in general will not give such a vertex algebra 
again.

In the present paper, we introduce a cohomology theory of grading-restricted
vertex algebras (including vertex operator algebras).
To overcome the difficulties in the proposal in \cite{KV}
mentioned above, our main new idea is 
to consider, instead of
linear maps from the tensor powers of the vertex algebra to 
a module for the algebra,  linear maps from the 
tensor powers of the vertex algebra
to suitable spaces of ``rational functions valued in 
the algebraic completion of the $V$-module'' such that they are ``composable''
with $m$ vertex operators in a natural sense and satisfy certain
other natural properties. 
These linear maps form a chain complex but it is still not a correct chain complex 
for the grading-restricted vertex algebra because the commutativity 
property for the vertex algebra has not been taken into consideration.
The correct chain complex for our cohomology 
is a subcomplex of this complex obtained by using shuffles 
in analogy with the construction of the Harrison cochain complex of
a commutative associative algebra from its Hochschild 
cochain complex. One subtle complication arising from the functions
mentioned above is that we have to  carefully
address the issue of convergence when we compose
these linear maps with vertex operators. In particular, for each $n\in \N$,
we have an inverse system $\{H^{n}_{m}(V, W)\}_{m\in \Z_{+}}$ of 
$n$-th cohomologies and an additional $n$-th cohomology $H_{\infty}^{n}(V, W)$
of a grading-restricted
vertex algebra $V$ with coefficients in a grading-restricted generalized 
$V$-module $W$ 
such that $H_{\infty}^{n}(V, W)$ is isomorphic to the inverse limit 
of the inverse system $\{H^{n}_{m}(V, W)\}_{m\in \Z_{+}}$. 
In the case of $n=2$, there is an additional 
second cohomology denoted by $H^{2}_{\frac{1}{2}}(V, W)$ which will be shown in
a sequel \cite{Hu3} to the present paper 
to correspond to what we call square-zero extensions of $V$ and to first order 
deformations of $V$ when $W=V$. 

The ideas and constructions in the present paper can also be applied to 
grading-restricted open-string vertex algebras (see \cite{HK1}) and grading-restricted
full field algebras (see \cite{HK2}) to introduce and study cohomologies
for these algebras. We shall present these cohomology theories 
in future publications. 

Note that for open-string vertex algebras and full field algebras, we have to 
work with complex variables, not formal variables. In particular, these algebras
are defined over only the field of complex numbers. In this paper,
we present our cohomology theory of
grading-restricted vertex algebras only over the field of complex numbers
so that it will be easy for us to
generalize the definitions and results given in the 
present paper to these algebras.
However, the cohomologies introduced in the present paper
can be defined and studied 
for grading-restricted vertex algebras over an arbitrary field $\mathbb{F}$
of characteristic $0$. In fact, to define and study
these cohomologies for such grading-restricted 
vertex algebras over $\mathbb{F}$, we need only replace rational functions 
with only possible poles at $z_{i}=z_{j}$ for $i\ne j$ by the localization of the 
polynomial ring $\mathbb{F}[z_{1}, \dots, z_{n}]$ by the first order polynomials 
$z_{i}-z_{j}$ for $i\ne j$ and replace series absolutely convergent to such rational 
functions in certain regions by series which are expansions of the elements of 
the localization corresponding to the expansions of such rational functions 
in the regions. See \cite{FHL} for discussions on formal rational functions and 
their expansions over such a field $\mathbb{F}$.

In a sequel \cite{Hu3} to the present paper, we shall 
show that for any $m\in \Z_{+}$, 
the first cohomology $H^{1}_{m}(V, W)$ of a grading-restricted
vertex algebra $V$ with coefficients in a grading-restricted
generalized $V$-module $W$ is linearly isomorphic 
to the space of derivations from $V$ to $W$.
We shall also
show that the 
second cohomology $H^{2}_{\frac{1}{2}}(V, W)$ of $V$ with coefficients in $W$ 
corresponds bijectively to the set of equivalence classes of square-zero extensions of 
$V$ by $W$ and the 
second cohomology $H^{2}_{\frac{1}{2}}(V, V)$ of $V$ with coefficients in $V$
corresponds bijectively to the set of equivalence classes of
first order deformations 
of $V$. 

At this moment, the author still does not 
have any vanishing theorem or duality theorem for the cohomologies introduced in the 
present paper. It is even not clear whether $H_{\infty}^{n}(V, W)$ vanishes when
$n$ is large. These are important research topics for the future development 
and applications of this cohomology theory.

This paper is organized as follows: In Section 2, we recall the notions of 
grading-restricted vertex algebra and grading-restricted generalized module
and also some useful results. In Section 3,
we introduce and study $\overline{W}$-valued rational functions for
a grading-restricted generalized module for a grading-restricted 
vertex algebra. These functions are 
crucial to our cohomology theory.
We present our cohomology theory in Section 4. 
%In Section 5, we discuss the long exact sequences
%of cohomologies induced from short exact sequences for grading-restricted
%vertex algebras or modules. 
%In Section 6, the bijection between
%the first cohomology of a grading-restricted 
%vertex algebra $V$ with coefficients in a $V$-module $W$ and 
%the space of derivations from $V$ to $W$ is given. 
%In Section 6,
%we show that the 
%second cohomology of $V$ with coefficients in $W$ 
%corresponds bijectively to the set of equivalence classes of square-zero extensions of 
%$V$ by $W$ and that the 
%second cohomology of $V$ with coefficients in $V$ 
%corresponds bijectively to the set of equivalence classes of
%first order deformations 
%of $V$.

The cohomologies introduced in the present paper 
was first presented in a talk by the author in 
Cao Xi-Hua Algebra Forum in East China Normal University
on June 1, 2010. 

\paragraph{Acknowledgments}
The author is grateful for partial 
support {}from NSF grant PHY-0901237.

\renewcommand{\theequation}{\thesection.\arabic{equation}}
\renewcommand{\thethm}{\thesection.\arabic{thm}}
\setcounter{equation}{0}
\setcounter{thm}{0}

\section{Grading-restricted 
vertex algebras and modules}

In this section, we give the definitions of grading-restricted 
vertex algebra and grading-restricted generalized module and 
discuss their basic properties. 
As is mentioned in the introduction,
we shall work  only over the field $\C$ of complex numbers in this paper.
In particular, all vector spaces are over $\C$.

A large part of the material in this 
section is from \cite{FHL} but we shall use the duality properties 
instead of the Jacobi identity in this paper. 
Below we recall the definition of 
grading-restricted vertex algebra using the duality properties as the main axiom. 

By a rational function of $z_{1}, \dots, z_{n}$, we mean 
a function of $z_{1}, \dots, z_{n}$ of the form 
$$f(z_{1}, \dots, z_{n})=\frac{P(z_{1}, \dots, z_{n})}{Q(z_{1}, \dots, z_{n})},$$
where $P(z_{1}, \dots, z_{n})$ and $Q(z_{1}, \dots, z_{n})$
are polynomials in $z_{1}, \dots, z_{n}$. 
If the polynomials $P(z_{1}, \dots, z_{n})$ and $Q(z_{1}, \dots, z_{n})$
have no common factors, then for a linear factor $g(z_{1}, \dots, z_{n})$ of 
$Q(z_{1}, \dots, z_{n})$, we say that {\it $f(z_{1}, \dots, z_{n})$
has poles at the set of zeros of $g(z_{1}, \dots, z_{n})$}
and the maximal power of $g(z_{1}, \dots, z_{n})$ in 
$Q(z_{1}, \dots, z_{n})$
is called the {\it order} of these poles. By a {\it rational function
with the only possible poles at a set of points in $\C^{n}$}, we mean 
a rational function of the form above such that 
$P(z_{1}, \dots, z_{n})$ and $Q(z_{1}, \dots, z_{n})$
have no common factors, $Q(z_{1}, \dots, z_{n})$ is a product of 
linear factors whose zeros are contained in that set of points in 
$\C^{n}$.

In the following definitions and in the rest of this paper,
$x, x_{1}, x_{2}, \dots$ are formal commuting variables 
and $z, z_{1}, z_{2}, \dots$ are complex numbers or complex 
variables.

\begin{defn}
{\rm A grading-restricted vertex algebra is 
a $\Z$-graded vector space $V=\coprod_{n\in \mathbb{Z}}V_{(n)}$
equipped with
a {\it vertex operator map} 
\begin{eqnarray*}
Y:V\otimes V&\to& V[[x, x^{-1}]], \nn
u\otimes v&\mapsto & Y_{V}(u, x)v=\sum_{n\in \Z}(Y_{V})_{n}(u)vx^{-n-1},
\end{eqnarray*}
a {\it vacuum} $\mathbf{1}\in V_{(0)}$
satisfying the following conditions:

\begin{enumerate}

\item {\it Grading restriction condition}:
For $n\in \Z$, $\dim V_{(n)}<\infty$ and when $n$ is sufficiently negative,
$V_{(n)}=0$. 

\item {\it Lower-truncation condition for vertex operators}:
For $u, v\in V$, $Y_{V}(u, x)v$ contain only finitely many negative 
power terms, that is, $Y_{V}(u, x)v\in V((x))$ (the space of formal 
Laurent series in $x$ with coefficients in $V$ and with finitely 
many negative power terms).  

\item {\it Identity property}: 
Let $1_{V}$ be the identity operator on $V$. Then 
$Y_{V}(\mathbf{1}, x)=1_{V}$.

\item {\it Creation property}: For $u\in V$, $Y_{V}(u, x)\mathbf{1}\in V[[x]]$
and $\lim_{x\to 0}Y_{V}(u, x)\mathbf{1}=u$.

\item {\it Duality}: For $u_{1}, u_{2}, v\in V$, 
$v'\in V'=\coprod_{n\in \mathbb{Z}}V_{(n)}^{*}$, the series 
\begin{eqnarray*}
&\langle v', Y_{V}(u_{1}, z_{1})Y_{V}(u_{2}, z_{2})v\rangle,&\\
&\langle v', Y_{V}(u_{2}, z_{2})Y_{V}(u_{1}, z_{1})v\rangle,&\\
&\langle v', Y_{V}(Y_{V}(u_{1}, z_{1}-z_{2})u_{2}, z_{2})v\rangle&
\end{eqnarray*}
are absolutely convergent
in the regions $|z_{1}|>|z_{2}|>0$, $|z_{2}|>|z_{1}|>0$,
$|z_{2}|>|z_{1}-z_{2}|>0$, respectively, to a common rational function 
in $z_{1}$ and $z_{2}$ with the only possible poles at $z_{1}, z_{2}=0$ and 
$z_{1}=z_{2}$. 

\item {\it $L(0)$-bracket formula}: Let $L_{V}(0): V\to V$ 
be defined by $L_{V}(0)v=nv$ for $v\in V_{(n)}$. Then
$$[L_{V}(0), Y_{V}(v, x)]=Y_{V}(L_{V}(0)v, x)+x\frac{d}{dx}Y_{V}(v, x)$$
for $v\in V$.

\item {\it $L(-1)$-derivative property}: 
Let $L_{V}(-1): V\to V$ be the operator 
given by 
$$L_{V}(-1)v=\res_{x}x^{-2}Y_{V}(v, x)\one=Y_{-2}(v)\one$$
for $v\in V$. Then for $v\in V$, 
$$\frac{d}{dx}Y_{V}(u, x)=Y_{V}(L_{V}(-1)u, x)=[L_{V}(-1), Y_{V}(u, x)].$$

\end{enumerate}}
\end{defn}

\begin{defn}
{\rm A {\it grading-restricted generalized $V$-module} is a vector space 
$W$ equipped with a vertex operator map 
\begin{eqnarray*}
Y_{W}: V\otimes W&\to& W[[x, x^{-1}]],\nn
u\otimes w&\mapsto & Y_{W}(u, x)w=\sum_{n\in \Z}(Y_{W})_{n}(u)wx^{-n-1}
\end{eqnarray*}
and linear operators $L_{W}(0)$ and $L_{W}(-1)$ on $W$ satisfying the following
conditions:

\begin{enumerate}

\item {\it Grading restriction condition}:
The vector space $W$ is $\C$-graded, that is, 
$W=\coprod_{n\in \mathbb{C}}W_{(n)}$, such that 
$W_{(n)}=0$ when the real part of $n$ is sufficiently negative. 

\item {\it Lower-truncation condition for vertex operators}:
For $u\in V$ and $w\in W$, $Y_{W}(u, x)w$ contain only finitely many negative 
power terms, that is, $Y_{W}(u, x)w\in W((x))$.

\item {\it Identity property}: 
Let $1_{W}$ be the identity operator on $W$. Then 
$Y_{W}(\mathbf{1}, x)=1_{W}$.

\item {\it Duality}: For $u_{1}, u_{2}\in V$, $w\in W$,
$w'\in W'=\coprod_{n\in \mathbb{Z}}W_{(n)}^{*}$, the series 
\begin{eqnarray*}
&\langle w', Y_{W}(u_{1}, z_{1})Y_{W}(u_{2}, z_{2})w\rangle,&\\
&\langle w', Y_{W}(u_{2}, z_{2})Y_{W}(u_{1}, z_{1})w\rangle,&\\
&\langle w', Y_{W}(Y_{V}(u_{1}, z_{1}-z_{2})u_{2}, z_{2})w\rangle&
\end{eqnarray*}
are absolutely convergent
in the regions $|z_{1}|>|z_{2}|>0$, $|z_{2}|>|z_{1}|>0$,
$|z_{2}|>|z_{1}-z_{2}|>0$, respectively, to a common rational function 
in $z_{1}$ and $z_{2}$ with the only possible poles at $z_{1}, z_{2}=0$ and 
$z_{1}=z_{2}$. 

\item {\it $L_{W}(0)$-bracket formula}: For  $v\in V$,
$$[L_{W}(0), Y_{W}(v, x)]=Y_{W}(L(0)v, x)+x\frac{d}{dx}Y_{W}(v, x).$$

\item {\it $L(0)$-grading property}: For $w\in W_{(n)}$, there exists
$N\in \Z_{+}$ such that $(L_{W}(0)-n)^{N}w=0$.

\item {\it $L(-1)$-derivative property}: For $v\in V$, 
$$\frac{d}{dx}Y_{W}(u, x)=Y_{W}(L_{V}(-1)u, x)=[L_{W}(-1), Y_{W}(u, x)].$$

\end{enumerate}
}
\end{defn}

Since in this paper, we shall always consider grading-restricted 
generalized $V$-modules,
for simplicity, we shall call them simply $V$-modules. 

If a meromorphic function $f(z_{1}, \dots, z_{n})$ on a region in $C^{n}$
can be analytically extended to a rational function in $z_{1}, \dots, z_{n}$, 
we shall
use $R(f(z_{1}, \dots, z_{n}))$ to denote this rational function. 

\begin{rema}
{\rm Let $V$ be a grading-restricted vertex algebra and $W$ a $V$-module. 
Then the 
duality axiom can be rewritten as:
For $u_{1}, u_{2}\in V$, $w\in W$, 
$w'\in W'$, 
\begin{eqnarray*}
R(\langle w', Y_{W}(u_{1}, z_{1})Y_{W}(u_{2}, z_{2})w\rangle)
&=&R(\langle w', Y_{W}(u_{2}, z_{2})Y_{W}(u_{1}, z_{1})w\rangle)\nn
&=&R(\langle w', Y_{W}(Y_{V}(u_{1}, z_{1}-z_{2})u_{2}, z_{2})w\rangle).
\end{eqnarray*}}
\end{rema}

The following result was proved in \cite{FHL} (Proposition 3.5.1 in \cite{FHL}):

\begin{prop}\label{n-comm}
For $v_{1}, \dots, v_{n}\in V$, $w\in W$ and
$w'\in W'$, 
$$\langle w', Y_{W}(v_{1}, z_{1})\cdots Y_{W}(v_{n}, z_{n})w\rangle$$
is absolutely convergent in the region $|z_{1}|>\cdots >|z_{n}|>0$ 
to a rational function 
$$R(\langle w', Y_{W}(v_{1}, z_{1})\cdots Y_{W}(v_{n}, z_{n})w\rangle)$$
in $z_{1}, \dots, z_{n}$ with the only possible poles at $z_{i}=z_{j}$, $i\ne j$,
and $z_{i}=0$. Moreover, the following commutativity holds: For 
$\sigma\in S_{n}$,
\begin{eqnarray*}
\lefteqn{R(\langle w', Y_{W}(v_{1}, z_{1})\cdots Y_{W}(v_{n}, z_{n})w\rangle)}\nn
&&=R(\langle w', Y_{W}(v_{\sigma(1)}, z_{\sigma(1)})\cdots 
Y_{W}(v_{\sigma(n)}, z_{\sigma(n)})w\rangle).
\end{eqnarray*}
\end{prop}

The following result, though not explicitly stated in \cite{FHL},
was implicitly given in Subsection 3.5 in \cite{FHL}: 

\begin{prop}\label{n-asso}
For $v_{1}, \dots, v_{n}\in V$, $w\in W$, 
$w'\in W'$ and $i=1, \dots, n-1$, 
\begin{eqnarray*}
\lefteqn{\langle w', Y_{W}(v_{1}, z_{1})\cdots 
Y_{W}(v_{i-1}, z_{i-1})\cdot}\nn
&&\quad\quad\quad \cdot 
Y_{W}(Y_{V}(v_{i}, z_{i}-z_{i+1})v_{i+1}, z_{i+1})
Y_{W}(v_{i+2}, z_{i+2})
\cdots Y_{W}(v_{n}, z_{n})w\rangle
\end{eqnarray*}
is absolutely convergent in the region given by $|z_{1}|>\cdots >|z_{i-1}|>|z_{i+1}|>
\cdots >|z_{n}|>0$, $|z_{i+1}|>|z_{i}-z_{i+1}|>0$ and 
$|z_{k}-z_{i+1}|>|z_{i}-z_{i+1}|>0$ for $k\ne i, i+1$
to a rational function 
\begin{eqnarray*}
\lefteqn{R(\langle w', Y_{W}(v_{1}, z_{1})\cdots 
Y_{W}(v_{i-1}, z_{i-1})\cdot}\nn
&&\quad\quad\quad \cdot Y_{W}(Y_{V}(v_{i}, z_{i}-z_{i+1})v_{i+1}, z_{i+1})
Y_{W}(v_{i+2}, z_{i+2})
\cdots Y_{W}(v_{n}, z_{n})w\rangle)
\end{eqnarray*}
with the only possible poles at $z_{i}=z_{j}$, $i\ne j$,
and $z_{i}=0$. Moreover, the following associativity holds: 
\begin{eqnarray*}
\lefteqn{R(\langle w', Y_{W}(v_{1}, z_{1})\cdots Y_{W}(v_{i-1}, z_{i-1})\cdot}\nn
&&\quad\quad\quad \cdot Y_{W}(v_{i}, z_{i})Y_{W}(v_{i+1}, z_{i+1})
Y_{W}(v_{i+2}, z_{i+2})
\cdots Y_{W}(v_{n}, z_{n})w\rangle)\nn
&&=R(\langle w', Y_{W}(v_{1}, z_{1})\cdots 
Y_{W}(v_{i-1}, z_{i-1})\cdot\nn
&&\quad\quad\quad \cdot Y_{W}(Y_{V}(v_{i}, z_{i}-z_{i+1})v_{i+1}, z_{i+1})
Y_{W}(v_{i+2}, z_{i+2})
\cdots Y_{W}(v_{n}, z_{n})w\rangle).
\end{eqnarray*}
\end{prop}

Recall from Subsection 5.6 in \cite{FHL} the linear map
\begin{eqnarray*}
Y_{WV}^{W}: W\otimes V&\to& W[[z, z^{-1}]]\nn
w\otimes v&\mapsto& Y_{WV}^{W}(w, z)v
\end{eqnarray*}
defined by
$$Y_{WV}^{W}(w, z)v=e^{zL(-1)}Y_{W}(v, -z)w$$
for $v\in V$ and $w\in W$. The following result is a special case
of Theorem 6.6.2 in \cite{FHL}:

\begin{prop}\label{w-n-comm}
For $v_{1}, \dots, v_{i-1}, v_{i+1}, \dots, v_{n}, v\in V$, $w\in W$ 
and $w'\in W'$,
$$\langle w', Y_{W}(v_{1}, z_{1})\cdots Y_{W}(v_{i-1}, z_{i-1})
Y_{WV}^{W}(w, z_{i}) Y_{V}(v_{i+1}, z_{i+1})\cdots Y_{V}(v_{n}, z_{n})v\rangle$$
is absolutely convergent in the region $|z_{1}|>\cdots >|z_{n}|>0$ 
to a rational function 
$$R(\langle w', Y_{W}(v_{1}, z_{1})\cdots Y_{W}(v_{i-1}, z_{i-1})
Y_{WV}^{W}(w, z_{i}) Y_{V}(v_{i+1}, z_{i+1})\cdots Y_{V}(v_{n}, z_{n})v\rangle)$$
in $z_{1}, \dots, z_{n}$
with the only possible poles at $z_{i}=z_{j}$, $i\ne j$,
and $z_{i}=0$. Moreover, the following commutativity holds: For 
$\sigma\in S_{n}$,
\begin{eqnarray*}
\lefteqn{R(\langle w', Y_{W}(v_{1}, z_{1})\cdots Y_{W}(v_{i-1}, z_{i-1})
Y_{WV}^{W}(w, z_{i}) Y_{V}(v_{i+1}, z_{i+1})\cdots Y_{V}(v_{n}, z_{n})v\rangle)}\nn
&&=R(\langle w', Y_{W}(u_{\sigma(1)}, z_{\sigma(1)})\cdots 
Y_{W}(u_{\sigma^{-1}(i)-1}, z_{\sigma^{-1}(i)-1})\nn
&&\quad\quad\quad\quad\quad\;\; \cdot
Y_{WV}^{W}(w, z_{i}) Y_{V}(u_{\sigma(\sigma^{-1}(i)+1)}, z_{\sigma^{-1}(i)+1})\cdots
Y_{W}(v_{\sigma(n)}, z_{\sigma(n)})v\rangle).
\end{eqnarray*}
\end{prop}

The following result, though not explicitly stated in \cite{FHL},
was implicitly given in Section 5.6 in \cite{FHL}: 

\begin{prop}\label{w-n-asso}
Let $v_{1}, \dots, v_{i-1}, v_{i+1}, \dots, v_{n}, v\in V$, $w\in W$ 
and $w'\in W'$. 
Let $u_{k}=v_{k}$ for $k=1, \dots, i-1, i+1, \dots, n$, $u_{i}=w$
and $1\le j\le n$. Let $\mathcal{Y}_{1}, \dots, \mathcal{Y}_{n}$ 
be $Y_{W}$, $Y_{WV}^{W}$ or $Y_{V}$ such that the expressions below are uniquely 
defined. Then  
\begin{eqnarray*}
\lefteqn{\langle w', \mathcal{Y}_{1}(u_{1}, z_{1})\cdots 
\mathcal{Y}_{j-1}(u_{j-1}, z_{j-1})\cdot}\nn
&&\quad \cdot 
\mathcal{Y}_{j}(\mathcal{Y}_{j+1}(u_{j}, z_{j}-z_{j+1})u_{j+1}, z_{j+1})
\mathcal{Y}_{j+2}(u_{j+2}, z_{j+2})
\cdots \mathcal{Y}_{n}(u_{n}, z_{n})v\rangle
\end{eqnarray*}
is
absolutely convergent in the region given by $|z_{1}|>\cdots >|z_{j-1}|>|z_{j+1}|>
\cdots >|z_{n}|>0$, $|z_{j+1}|>|z_{j}-z_{j+1}|>0$ and 
$|z_{k}-z_{j+1}|>|z_{j}-z_{j+1}|>0$ for $k\ne j, j+1$
to a rational function 
\begin{eqnarray*}
\lefteqn{R(\langle w', \mathcal{Y}_{1}(u_{1}, z_{1})\cdots 
\mathcal{Y}_{j-1}(u_{j-1}, z_{j-1})\cdot}\nn
&&\quad \cdot 
\mathcal{Y}_{j}(\mathcal{Y}_{j+1}(u_{j}, z_{j}-z_{j+1})u_{j+1}, z_{j+1})
\mathcal{Y}_{j+2}(u_{j+2}, z_{j+2})
\cdots \mathcal{Y}_{n}(u_{n}, z_{n})v\rangle)
\end{eqnarray*}
in $z_{1}, \dots, z_{n}$
with the only possible poles at $z_{i}=z_{j}$, $i\ne j$,
and $z_{i}=0$. Moreover, the following associativity holds: 
For $j=1, \dots, i-2$ or $j=i+1, \dots, n-1$, 
\begin{eqnarray*}
\lefteqn{R(\langle w', Y_{W}(v_{1}, z_{1})\cdots Y_{W}(v_{i-1}, z_{i-1})\cdot}\nn
&&\quad \cdot 
Y_{WV}^{W}(w, z_{i}) Y_{V}(v_{i+1}, z_{i+1})\cdots Y_{V}(v_{n}, z_{n})v\rangle)\nn
&&=R(\langle w', \mathcal{Y}_{1}(u_{1}, z_{1})\cdots 
\mathcal{Y}_{j-1}(u_{j-1}, z_{j-1})\cdot\nn
&&\quad \cdot 
\mathcal{Y}_{j}(\mathcal{Y}_{j+1}(u_{j}, z_{j}-z_{j+1})u_{j+1}, z_{j+1})
\mathcal{Y}_{j+2}(u_{j+2}, z_{j+2})
\cdots \mathcal{Y}_{n}(u_{n}, z_{n})v\rangle).
\end{eqnarray*}
\end{prop}

Let $\overline{W}$ be the algebraic completion of $W$, that is, 
$\overline{W}=\prod_{n\in \C}W_{(n)}=(W')^{*}$. For $n\in \Z_{+}$,
let $F_{n}\C$ be the configuration space of $n$ points in $\C$, that is,
$$F_{n}\C=\{(z_{1}, \dots, z_{n})\in \C^{n}\;|\; z_{i}\ne z_{j}, i\ne j\}.$$
For each $(z_{1}, \dots, z_{n}, \zeta)\in
F_{n+1}\C$, $v_{1}, \dots, v_{n}\in V$, $w\in W$ and $w'\in W'$, we have an element 
$$E(Y_{W}(v_{1}, z_{1})\cdots Y_{W}(v_{n}, z_{n})Y_{WV}^{W}(w, \zeta)\one)\in \overline{W}$$
given by 
\begin{eqnarray*}
\lefteqn{\langle w',E(Y_{W}(v_{1}, z_{1})\cdots Y_{W}(v_{n}, z_{n})Y_{WV}^{W}(w, \zeta)
\one)\rangle}\nn
&&=R(\langle w', Y_{W}(v_{1}, z_{1})\cdots Y_{W}(v_{n}, z_{n})Y_{WV}^{W}(w, \zeta)
\one\rangle).
\end{eqnarray*}
For $(z_{1}, \dots, z_{n}, \zeta)\in
F_{n+1}\C$, $v_{1}, \dots, v_{n}\in V$ and $w\in W$, 
set 
\begin{eqnarray*}
\lefteqn{(E^{(n, 1)}_{W}(v_{1}\otimes \cdots\otimes v_{n}; w))(z_{1}, \dots, z_{n}, \zeta)}\nn
&&=E(Y_{W}(v_{1}, z_{1})\cdots Y_{W}(v_{n}, z_{n})Y_{WV}^{W}(w, \zeta)
\one)\in \overline{W}.
\end{eqnarray*}
We also define
$$E^{(n)}_{W}(v_{1}\otimes \cdots\otimes v_{n}; w):
\{(z_{1}, \dots, z_{n})\in F_{n}\C\;|\; z_{i}\ne 0,\;i=1, \dots, n\}\to 
\overline{W}$$
by 
\begin{eqnarray*}
\lefteqn{(E^{(n)}_{W}(v_{1}\otimes \cdots\otimes v_{n}; w))
(z_{1}, \dots, z_{n})}\nn
&&=(E^{(n, 1)}_{W}(v_{1}\otimes \cdots\otimes v_{n}; w))(z_{1}, \dots, z_{n}, 0)\nn
&&=E(Y_{W}(v_{1}, z_{1})\cdots Y_{W}(v_{n}, z_{n})w).
\end{eqnarray*}

The next result was in fact proved in \cite{H1} and \cite{H2} but was formulated
using the geometry of spheres with punctures and local coordinates
in a more general setting. Here we formulate it without using the language
of geometry or operads and give a direct
proof. Given a $V$-module $W=\coprod_{n\in \C}W_{(n)}$, let $P_{n}: W\to W_{(n)}$
for $n\in \C$ be the projection from $W$ to $W_{(n)}$.

\begin{prop}\label{correl-fn}
For $k, l_{1}, \dots, l_{n+1} \in \Z_+$ and 
$v_{1}^{(1)},\dots, v_{l_{1}}^{(1)},\dots, v_{1}^{(n+1)}, \dots$,
$v_{l_{n+1}}^{(n+1)} \in V$, $w\in W$ and $w'\in W'$, the series 
\begin{eqnarray}\label{va-conv-axiom}
\lefteqn{\sum_{r_{1}, \ldots, r_{n}\in \Z, r_{n+1}\in \C} 
\langle w', (E^{(n, 1)}_{W}(P_{r_{1}} 
((E^{(l_{1})}_{V}(v_{1}^{(1)}\otimes
\cdots\otimes v_{l_{1}}^{(1)}; \one))(z_{1}^{(1)}, \dots, 
z_{l_{1}}^{(1)}))}  \nn
&& \qquad \otimes \cdots\otimes 
P_{r_n} ((E^{(l_{n})}_{V}(v_{1}^{(n)}\otimes  \cdots \otimes  v_{l_{n}}^{(n)};\one)) 
(z_{1}^{(n)}, \dots, z_{l_{n}}^{(n)}));\nn
&& \qquad P_{r_{n+1}}((E^{(l_{n+1})}_{W}(v_{1}^{(n+1)}
\otimes  \cdots \otimes  v_{l_{n+1}}^{(n+1)}; 
w))(z_{1}^{(n+1)}, \dots, z_{l_{n+1}}^{(n+1)}))))\nn
&& \qquad\qquad\qquad\qquad\qquad\qquad\qquad\qquad\qquad 
\qquad\qquad\quad(z^{(0)}_{1}, \dots, 
z^{(0)}_{n+1})\rangle\nn
\end{eqnarray}
converges absolutely to  
\begin{eqnarray}\label{sum}
\lefteqn{\langle w', (E^{(n)}_{W}(v_{1}^{(1)}\otimes  \cdots \otimes v_{l_{n+1}}^{(n+1)}; w))
(z_{1}^{(1)}+z^{(0)}_{1}, \dots,
z_{l_{1}}^{(1)}+z^{(0)}_{1}, }  \nn
&&\quad \quad  \quad  \quad  \quad  \quad  \quad  \quad  \quad  
\dots, 
z_{1}^{(n+1)}+z^{(0)}_{n+1},   \dots,
z_{l_{n+1}}^{(n+1)}+z^{(0)}_{n+1})\rangle. 
\end{eqnarray}
when $0<|z_p^{(i)}| + |z_q^{(j)}|< |z^{(0)}_i
-z^{(0)}_j|$ for $i,j=1, \dots, n+1$, $i\ne j$, $p=1, 
\dots,  l_i$, $q=1, \dots, l_j$.
\end{prop}
\pf
By definition,
$$\langle w', (E^{(n)}_{W}(v_{1}^{(1)}\otimes  \cdots \otimes v_{l_{n+1}}^{(n+1)}; w))
(z_{1}, \dots, z_{l_{1}+\cdots +l_{n+1}})\rangle$$
is a rational function in $z_{1}, \dots, z_{l_{1}+\cdots +l_{n+1}}$ with the 
only possible poles at $z_{i}=0$ or $z_{i}=z_{j}$, $i\ne j$.
Thus for fixed $z_p^{(i)}\in \C^{\times}$, $p=1, 
\dots,  l_i$, $i=1, \dots, n+1$ and $z_{i}^{(0)}$ for $i=1, \dots, n+1$
satisfying $0<|z_p^{(i)}| + |z_q^{(j)}|< |z^{(0)}_i
-z^{(0)}_j|$ for $i,j=1, \dots, n+1$, $i\ne j$, $p=1, 
\dots,  l_i$, $q=1, \dots, l_j$, the function 
\begin{eqnarray}\label{sum-t}
\lefteqn{\langle w', (E^{(n)}_{W}(t_{1}^{L(0)}v_{1}^{(1)}\otimes \cdots
\otimes t_{1}^{L(0)}v_{l_{1}}^{(1)}\otimes  \cdots\otimes 
t_{n+1}^{L(0)}v_{1}^{(n+1)}\otimes 
\cdots \otimes t_{n+1}^{L(0)}v_{l_{n+1}}^{(n+1)}; w))}
\nn
&&\quad  \quad\quad\quad\quad (t_{1}z_{1}^{(1)}+z^{(0)}_{1}, \dots,
t_{1}z_{l_{1}}^{(1)}+z^{(0)}_{1},  \nn
&& \quad  \quad\quad\quad  \quad\quad\quad  \quad\quad\quad\quad
\dots, 
t_{n+1}z_{1}^{(n+1)}+z^{(0)}_{n+1},   \dots,
t_{n+1}z_{l_n}^{(n+1)}+z^{(0)}_{n+1})\rangle. \nn
\end{eqnarray}
of $t_{1}, \dots, t_{n+1}$ has an expansion
as a Laurent series in $t_{1}, \dots, t_{n+1}$ when $(t_{1}, \dots, t_{n+1})$ 
are in the direct product of some annuli containing $1$. 
Using induction and the associativity for $V$ and $W$ repeatedly, we
see that the coefficients of this Laurent expansion are the same 
as the coefficients of the formal Laurent series
\begin{eqnarray}\label{va-conv-axiom-t}
\lefteqn{\sum_{r_{1}, \ldots, r_{n+1}\in \Z} \langle w', (E^{(n, 1)}_{W}(P_{r_{1}} 
((E^{(l_{1})}_{V}(v_{1}^{(1)}\otimes
\cdots\otimes v_{l_{1}}^{(1)}; \one))(z_{1}^{(1)}, \dots, 
z_{l_{1}}^{(1)}))}  \nn
&& \otimes \cdots\otimes 
P_{r_n} ((E^{(l_{n})}_{V}(v_{1}^{(n)}\otimes  \cdots \otimes  v_{l_{n}}^{(n)};\one)) 
(z_{1}^{(n)}, \dots, z_{l_{n}}^{(n)})); \nn
&& P_{r_{n+1}}((E^{(l_{n+1})}_{W}(v_{1}^{(n+1)}
\otimes  \cdots \otimes  v_{l_{n+1}}^{(n+1)}; 
w))(z_{1}^{(n+1)}, \dots, z_{l_{n+1}}^{(n+1)}))))\nn
&& \qquad\qquad\qquad\qquad\qquad\qquad\qquad\qquad\qquad  (z^{(0)}_{1}, \dots, 
z^{(0)}_{n+1})\rangle t_{1}^{r_{1}}\cdots t_{n}^{r_{n}}.\nn
&&
\end{eqnarray}
Thus (\ref{va-conv-axiom-t}) is absolutely convergent to (\ref{sum-t})
in the region where (\ref{sum-t}) has a Laurent expansion. In particular,
when $t_{1}=\cdots=t_{n+1}=1$, we obtain that (\ref{va-conv-axiom}) 
is absolutely convergent to (\ref{sum}).
\epfv

For each $(\zeta, z_{1}, \dots, z_{n})\in
F_{n+1}\C$, $v_{1}, \dots, v_{n}\in V$, $w\in W$ and $w'\in W'$, we have an element 
$$E(Y_{WV}^{W}(w, \zeta)Y_{V}(v_{1}, z_{1})\cdots Y_{V}(v_{n}, z_{n})\one)\in \overline{W}$$
given by 
\begin{eqnarray*}
\lefteqn{\langle w',E(Y_{WV}^{W}(w, \zeta)Y_{V}(v_{1}, z_{1})\cdots 
Y_{V}(v_{n}, z_{n})\one)\rangle}\nn
&&=R(\langle w', Y_{WV}^{W}(w, \zeta)Y_{V}(v_{1}, z_{1})\cdots 
Y_{V}(v_{n}, z_{n})\one\rangle).
\end{eqnarray*}
For $(\zeta, z_{1}, \dots, z_{n})\in
F_{n+1}\C$, $v_{1}, \dots, v_{n}\in V$ and $w\in W$, 
set 
\begin{eqnarray*}
\lefteqn{(E^{W; (1, n)}_{WV}(w; v_{1}\otimes \cdots\otimes v_{n}))(\zeta, 
z_{1}, \dots, z_{n})}\nn
&&=E(Y_{WV}^{W}(w, \zeta)Y_{V}(v_{1}, z_{1})\cdots 
Y_{V}(v_{n}, z_{n})\one)\in \overline{W}.
\end{eqnarray*}
We have:

\begin{prop}\label{wv<=>vw}
For $(\zeta, z_{1}, \dots, z_{n})\in
F_{n+1}\C$, $v_{1}, \dots, v_{n}\in V$ and $w\in W$,
\begin{eqnarray*}
\lefteqn{(E^{W; (1, n)}_{WV}(w; v_{1}\otimes \cdots\otimes v_{n}))
(\zeta, z_{1}, \dots, z_{n})}\nn
&&=(E^{(n, 1)}_{W}(v_{1}\otimes \cdots\otimes v_{n}; w))(z_{1}, \dots, z_{n}, \zeta).
\end{eqnarray*}
\end{prop}
\pf
For $(\zeta, z_{1}, \dots, z_{n})\in
F_{n+1}\C$, $v_{1}, \dots, v_{n}\in V$, $w\in W$ and $w'\in W'$, 
\begin{eqnarray*}
\lefteqn{\langle w', (E^{W; (1, n)}_{WV}(w; v_{1}\otimes \cdots\otimes v_{n}))
(\zeta, z_{1}, \dots, z_{n})\rangle}\nn
&&=R(\langle w', Y^{W}_{WV}(w, \zeta)Y_{V}(v_{1},  z_{1})\cdots
Y_{V}(v_{n}, z_{n})\one\rangle)\nn
&&=R(\langle w', e^{\zeta L(-1)}Y_{W}(Y_{V}(v_{1}, z_{1})\cdots
Y_{V}(v_{n}, z_{n})\one, -\zeta)w\rangle)\nn
&&=R(\langle w', e^{\zeta L(-1)}Y_{W}(v_{1}, z_{1}-\zeta)\cdots
Y_{W}(v_{n}, z_{n}-\zeta)w\rangle)\nn
&&=R(\langle w', Y_{W}(v_{1}, z_{1})\cdots
Y_{W}(v_{n}, z_{n})e^{\zeta L(-1)}w\rangle)\nn
&&=R(\langle w', Y_{W}(v_{1}, z_{1})\cdots
Y_{W}(v_{n}, z_{n})Y_{WV}^{W}(w, \zeta)\one\rangle)\nn
&&=\langle w', 
(E^{(n, 1)}_{W}(v_{1}\otimes \cdots\otimes v_{n}; w))(z_{1}, \dots, z_{n}, \zeta)\rangle.
\end{eqnarray*}
\epfv

We also define
$$E^{W; (n)}_{WV}(w; v_{1}\otimes \cdots\otimes v_{n}):
\{(z_{1}, \dots, z_{n})\in F_{n}\C\;|\; z_{i}\ne 0,\;i=1, \dots, n\}
\to \overline{W}$$
by 
$$(E^{W; (n)}_{WV}(w; v_{1}\otimes \cdots\otimes v_{n}))(z_{1}, \dots, z_{n})
=(E^{W; (1, n)}_{WV}(w; v_{1}\otimes \cdots\otimes v_{n}))(0, z_{1}, \dots, z_{n}).$$
Then by Proposition \ref{wv<=>vw}, 
$$E^{W; (n)}_{WV}(w; v_{1}\otimes \cdots\otimes v_{n})
=E^{(n)}_{W}(v_{1}\otimes \cdots\otimes v_{n}; w)$$
for $v_{1}, \dots, v_{n}\in V$ and $w\in W$.

\renewcommand{\theequation}{\thesection.\arabic{equation}}
\renewcommand{\thethm}{\thesection.\arabic{thm}}
\setcounter{equation}{0}
\setcounter{thm}{0}

\section{$\overline{W}$-valued rational functions}

Let $V$ be a grading-restricted 
vertex algebra and $W$ a $V$-module (recall our convention that 
a $V$-module means a grading-restricted generalized $V$-module 
in this paper). Recall the configuration spaces 
$$F_{n}\C=\{(z_{1}, \dots, z_{n})\in \C^{n}\;|\; z_{i}\ne z_{j}, i\ne j\}$$
for $n\in \Z_{+}$.

\begin{defn}
{\rm A {\it $\overline{W}$-valued rational function in $z_{1}, \dots, z_{n}$
with the only possible poles at 
$z_{i}=z_{j}$, $i\ne j$}
is a map 
\begin{eqnarray*}
f: \quad\quad\quad F_{n}\C
&\to& \overline{W}\nn
(z_{1}, \dots, z_{n})&\mapsto& f(z_{1}, \dots, z_{n})
\end{eqnarray*}
such that  for any $w'\in W'$, 
$$\langle w', f(z_{1}, \dots, z_{n})\rangle$$
is a rational function in $z_{1}, \dots, z_{n}$ 
with the only possible poles  at 
$z_{i}=z_{j}$, $i\ne j$.}
\end{defn} 

For simplicity, we shall call the map that we just defined  a {\it 
$\overline{W}$-valued rational function in
$z_{1}, \dots, z_{n}$} unless there might be other poles. 
Denote the space of all $\overline{W}$-valued rational functions in
$z_{1}, \dots, z_{n}$ by $\widetilde{W}_{z_{1}, \dots, z_{n}}$. 
We define a left action of $S_{n}$ on $\widetilde{W}_{z_{1}, \dots, z_{n}}$
by
$$(\sigma(f))(z_{1}, \dots, z_{n})=f(z_{\sigma(1)}, \dots, z_{\sigma(n)})$$
for $f\in \widetilde{W}_{z_{1}, \dots, z_{n}}$.

\begin{expl}
{\rm For $w\in W$, the $\overline{W}$-valued function
$E^{(n)}_{W}(v_{1}\otimes \cdots\otimes v_{n}; w)$
given by 
$$(E^{(n)}_{W}(v_{1}\otimes \cdots\otimes v_{n}; w))(z_{1}, \dots, z_{n})
=E(Y_{W}(v_{1}, z_{1})\cdots Y_{W}(v_{n}, z_{n})w)$$
in general might not be an element of $\widetilde{W}_{z_{1}, \dots, z_{n}}$, 
since there might be sigularities at $z_{i}=0$.  But for $w\in W$ 
such that $Y_{W}(v, x)w\in W[[x]]$, $E^{(n)}_{W}(v_{1}\otimes \cdots\otimes v_{n}; w)$ 
is indeed an element of $\widetilde{W}_{z_{1}, \dots, z_{n}}$. Since 
$$E^{W; (n)}_{WV}(w; v_{1}\otimes \cdots\otimes v_{n})
=E^{(n)}_{W}(v_{1}\otimes \cdots\otimes v_{n}; w),$$
$E^{W; (n)}_{WV}(w; v_{1}\otimes \cdots\otimes v_{n})$ is also 
an element of $\widetilde{W}_{z_{1}, \dots, z_{n}}$.
In particular, when $W=V$ and $w=\one$, 
$E^{(n)}_{V}(v_{1}\otimes \cdots\otimes v_{n}; \one)\in 
\widetilde{V}_{z_{1}, \dots, z_{n}}$.}
\end{expl}

For $z\in \C^{\times}$, we shall use $\log z$ to denote 
$\log |z|+i\arg z$, $0\le \arg z<2 \pi$. Let $(L_{W}(0))_{s}$ 
be the semisimple part of $L_{W}(0)$, that is, $(L_{W}(0))_{s}w=nw$ for $w\in W_{(n)}$. 
Since $W$ is a (grading-restricted generalized) $V$-module,
for any $z\in \C^{\times}$, 
\begin{eqnarray*}
z^{L_{W}(0)}&=&e^{(\log z)L_{W}(0)}\nn
&=&e^{(\log z)(L_{W}(0))_{s}}e^{(\log z)(L_{W}(0)-(L_{W}(0))_{s})}
\end{eqnarray*}
is a well-defined linear operator
on $\overline{W}$.

\begin{defn}
{\rm For $n\in \Z_{+}$, a linear map $\Phi: V^{\otimes n}\to 
\widetilde{W}_{z_{1}, \dots, z_{n}}$ is said to have
the {\it $L(-1)$-derivative property} if (i) 
\begin{eqnarray*}
\lefteqn{\frac{\partial}{\partial z_{i}}\langle w', 
(\Phi(v_{1}\otimes \cdots  \otimes v_{n}))(z_{1}, 
\dots, z_{n})\rangle}\nn
&&=\langle w', 
(\Phi(v_{1}\otimes \cdots \otimes v_{i-1}\otimes L_{V}(-1)v_{i}
\otimes v_{i+1}\otimes \cdots \otimes v_{n}))(z_{1}, \dots, z_{n})\rangle
\end{eqnarray*}
for $i=1, \dots, n$, $v_{1}, \dots, v_{n}\in V$ and
$w'\in W'$ and (ii) 
\begin{eqnarray*}
\lefteqn{\left(\frac{\partial}{\partial z_{1}}+\cdots +\frac{\partial}{\partial z_{n}}\right)
\langle w', 
(\Phi(v_{1}\otimes \cdots \otimes v_{n}))(z_{1}, \dots, z_{n})\rangle}\nn
&&\quad=\langle w', L_{W}(-1)
(\Phi(v_{1}\otimes \cdots  \otimes v_{n}))(z_{1}, 
\dots, z_{n})\rangle
\end{eqnarray*}
and  $v_{1}, \dots, v_{n}\in V$,
$w'\in W'$.
A linear map $\Phi: V^{\otimes n}\to 
\widetilde{W}_{z_{1}, \dots, z_{n}}$ is said to have
the {\it $L(0)$-conjugation property} if for $v_{1}, \dots, v_{n}\in V$,
$w'\in W'$, $(z_{1}, \dots, z_{n})\in F_{n}\C$ and $z\in \C^{\times}$ so that 
$(zz_{1}, \dots, zz_{n})\in F_{n}\C$,
\begin{eqnarray*}
\lefteqn{\langle w', z^{L_{W}(0)}
(\Phi(v_{1}\otimes \cdots \otimes v_{n}))(z_{1}, \dots, z_{n})\rangle}\nn
&&=\langle w', 
(\Phi(z^{L(0)}v_{1}\otimes \cdots \otimes  z^{L(0)}v_{n}))(zz_{1}, 
\dots, zz_{n})\rangle.
\end{eqnarray*}} 
\end{defn}

Note that since $L_{W}(-1)$ is a weight-one operator on $W$, 
for any $z\in \C$,  $e^{zL_{W}(-1)}$ is a well-defined linear operator
on $\overline{W}$. 

\begin{prop}
Let $\Phi: V^{\otimes n}\to 
\widetilde{W}_{z_{1}, \dots, z_{n}}$ be a linear map having
the $L(-1)$-derivative property. Then for $v_{1}, \dots, v_{n}\in V$,
$w'\in W'$, $(z_{1}, \dots, z_{n})\in F_{n}\C$, $z\in \C$ such that
$(z_{1}+z, \dots,  z_{n}+z)\in F_{n}\C$,
\begin{eqnarray*}
\lefteqn{\langle w', e^{zL_{W}(-1)}
(\Phi(v_{1}\otimes \cdots \otimes v_{n}))(z_{1}, \dots, z_{n})\rangle}\nn
&&=\langle w', 
(\Phi(v_{1}\otimes \cdots  \otimes v_{n}))(z_{1}+z, 
\dots, z_{n}+z)\rangle
\end{eqnarray*}
and for $v_{1}, \dots, v_{n}\in V$,
$w'\in W'$, $(z_{1}, \dots, z_{n})\in F_{n}\C$, $z\in \C$ and $1\le i\le n$ such that
$$(z_{1}, \dots, z_{i-1}, z_{i}+z, z_{i+1}, \dots, z_{n})\in F_{n}\C,$$
the power series expansion of 
\begin{equation}\label{expansion-fn}
\langle w', 
(\Phi(v_{1}\otimes \cdots  \otimes v_{n}))(z_{1}, 
\dots, z_{i-1}, z_{i}+z, z_{i+1}, \dots, z_{n})\rangle
\end{equation}
in $z$ is equal to the power series
\begin{equation}\label{power-series}
\langle w', 
(\Phi(v_{1}\otimes \cdots \otimes v_{i-1}\otimes e^{zL(-1)}v_{i}
\otimes v_{i+1}\otimes \cdots \otimes v_{n}))(z_{1}, \dots, z_{n})\rangle
\end{equation}
in $z$.
In particular, the power series (\ref{power-series}) in $z$ is absolutely convergent
to (\ref{expansion-fn}) in the disk $|z|<\min_{i\ne j}\{|z_{i}-z_{j}|\}$. 
\end{prop}
\pf
This result follows immediately from the definition of $L(-1)$-derivative property and
Taylor's theorem on power series expansions
of analytic functions.
\epfv 

We would like to take linear maps from tensor powers of $V$ to 
$\widetilde{W}_{z_{1}, \dots, z_{n}}$ to be cochains in our cohomology 
theory. But to define the coboundary operator, we have to compose cochains 
with vertex operators. However,
the images of vertex operator maps in general are not 
in the algebras or in the modules. They are in the algebraic completions
of the algebras or modules. This is one of the most subtle features 
of the theory of grading-restricted vertex algebras 
or vertex operator algebras. Because of this subtlety, we cannot 
compose vertex operators directly. Instead, we first write a series by projecting
an element of the algebraic completion of an algebra or a module 
to its homogeneous components, composing these homogeneous components
with other vertex operators 
and then taking the formal sum. If this formal sum
is absolutely convergent, then these operators
can be composed and we shall use the usual notation to denote the composition 
obtained from the sums of these series. See \cite{H2} for detailed discussions
in the case of vertex operator algebras. 

Since $\overline{W}$-valued rational functions above are valued in $\overline{W}$,
not $W$, and for $z\in \C^{\times}$, $u, v\in V$, $w\in W$, $Y_{V}(u, z)v\in \overline{V}$
and $Y_{W}(u, z)v\in \overline{W}$, in general, 
we might not be able to compose a linear map from 
a tensor power of $V$ to $\widetilde{W}_{z_{1}, \dots, z_{n}}$ with 
vertex operators.  So we have to consider linear maps
from tensor powers of $V$ to $\widetilde{W}_{z_{1}, \dots, z_{n}}$ 
such that these maps can be composed with vertex operators in the sense mentioned above. 

For a $V$-module $W=\coprod_{n\in \C}W_{(n)}$ and $m\in \C$, 
let $P_{m}: \overline{W}\to W_{(m)}$ be 
the projection from $\overline{W}$ to $W_{(m)}$. 

\begin{defn}\label{composable}
{\rm 
Let $\Phi: V^{\otimes n}\to 
\widetilde{W}_{z_{1}, \dots, z_{n}}$ be a linear map. For $m\in \N$, 
$\Phi$ is said to be {\it composable with $m$ vertex operators} if 
the following 
conditions are satisfied:

\begin{enumerate}

\item Let $l_{1}, \dots, l_{n}\in \Z_+$ such that $l_{1}+\cdots +l_{n}=m+n$,
$v_{1}, \dots, v_{m+n}\in V$ and $w'\in W'$. Set 
\begin{eqnarray}\label{psi-i}
\Psi_{i}
&=&(E^{(l_{i})}_{V}(v_{l_{1}+\cdots +l_{i-1}+1}\otimes
\cdots\otimes v_{l_{1}+\cdots +l_{i-1}+l_{i}}; \one))\nn
&&\quad\quad\quad\quad\quad\quad\quad\quad(z_{l_{1}+\cdots +l_{i-1}+1}-\zeta_{i}, \dots, 
z_{l_{1}+\cdots +l_{i-1}+l_{i}}-\zeta_{i})\nn
\end{eqnarray}
for $i=1, \dots, n$. Then there exist positive integers $N(v_{i}, v_{j})$
depending only on $v_{i}$ and $v_{j}$ for $i, j=1, \dots, k$, $i\ne j$ such that the series 
$$\sum_{r_{1}, \dots, r_{n}\in \Z}\langle w', 
(\Phi(P_{r_{1}}\Psi_{1}\otimes \dots\otimes
P_{r_n} \Psi_{n}))(\zeta_{1}, \dots, 
\zeta_{n})\rangle,$$
is absolutely convergent  when 
$$|z_{l_{1}+\cdots +l_{i-1}+p}-\zeta_{i}| 
+ |z_{l_{1}+\cdots +l_{j-1}+q}-\zeta_{i}|< |\zeta_{i}
-\zeta_{j}|$$
for $i,j=1, \dots, k$, $i\ne j$ and for $p=1, 
\dots,  l_i$ and $q=1, \dots, l_j$.
and the sum can be analytically extended to a
rational function
in $z_{1}, \dots, z_{m+n}$, independent of $\zeta_{1}, \dots, \zeta_{n}$,
with the only possible poles at 
$z_{i}=z_{j}$ of order less than or equal to 
$N(v_{i}, v_{j})$ for $i,j=1, \dots, k$,  $i\ne j$. 

\item For $v_{1}, \dots, v_{m+n}\in V$,  there exist 
positive integers $N(v_{i}, v_{j})$ depending only on $v_{i}$ and 
$v_{j}$ for $i, j=1, \dots, k$, $i\ne j$ such that for $w'\in W'$,
\begin{eqnarray*}
\lefteqn{\sum_{q\in \C}\langle w', 
(E^{(m)}_{W}(v_{1}\otimes \cdots\otimes v_{m}; }\nn
&&\quad\quad\quad
P_{q}((\Phi(v_{m+1}\otimes \cdots\otimes v_{m+n}))(z_{m+1}, \dots, z_{m+n})))
(z_{1}, \dots, z_{m})\rangle
\end{eqnarray*}
is absolutely convergent when $z_{i}\ne z_{j}$, $i\ne j$
$|z_{i}|>|z_{k}|>0$ for $i=1, \dots, m$ and 
$k=m+1, \dots, m+n$ and the sum can be analytically extended to a
rational function 
in $z_{1}, \dots, 
z_{m+n}$ with the only possible poles at 
$z_{i}=z_{j}$ of orders less than or equal to 
$N(v_{i}, v_{j})$ for $i, j=1, \dots, k$, $i\ne j$,. 

\end{enumerate}}
\end{defn}

\begin{rema}
{\rm In the first version of the present paper, we did not assume the existence of 
the positive integers $N(v_{i}, v_{j})$ for $i, j=1, \dots, k$, $i\ne j$. 
We did get a cohomology theory without such an assumption. On the other hand, since 
the correlation functions for grading-restricted vertex algebras do have this 
property, here we add this assumption. (In fact, the existence of $N(v_{i}, v_{i+1})$
can be seen immediately from Proposition \ref{n-asso} and the fact that 
$Y_{V}(v_{i}, z_{i}-z_{i+1})v_{i+1}$ contains only finitely many negative power terms
in $z_{i}-z_{i+1}$ (the lower-truncation condition). The existence of 
$N(v_{i}, v_{j})$ then follows from the existence of $N(v_{i}, v_{i+1})$
and Proposition \ref{n-comm}.) But we remark that the cohomology theory without 
this assumption might still be important in the future studies. We might 
call the cohomology theory without this assumption 
of the paper the {\it cohomology theory 
without upper bounds on orders of poles}.}
\end{rema}

We shall denote the rational functions in Conditions
1 and 2 of Definition \ref{composable}
by 
\begin{eqnarray*}
\lefteqn{R(\langle w', \Phi(E^{(l_{1})}_{V}(v_{1}\otimes \cdots 
\otimes v_{l_{1}}; \one)\otimes}\nn
&&\quad\quad\quad \cdots\otimes
E^{(l_{n})}_{V}(v_{l_{1}+\cdots +l_{n-1}+1}\otimes \cdots 
\otimes v_{l_{1}+\cdots +l_{n-1}+l_{n}}; \one))
(z_{1}, \dots, 
z_{m+n})\rangle)
\end{eqnarray*}
and 
$$R(\langle w', 
(E^{(m)}_{W}(v_{1}\otimes \cdots\otimes v_{m}; 
\Phi(v_{m+1}\otimes \cdots\otimes v_{m+n}))
(z_{1}, \dots, z_{m+n})\rangle),$$
respectively.

\begin{expl}\label{correl-fns}
{\rm For $w\in W$ satisfying $Y_{W}(v, x)w\in W[[x]]$, 
the $\overline{W}$-valued rational function 
$E^{(n)}_{W}(v_{1}\otimes \cdots\otimes v_{n}; w)$
for $v_{1}, \dots,  v_{n}\in V$
give a linear map 
\begin{eqnarray*}
E^{(n)}_{W;\;w}: \qquad \qquad V^{\otimes n}&\to &
\widetilde{W}_{z_{1}, \dots, z_{n}}\nn
v_{1}\otimes \cdots\otimes v_{n}& \mapsto 
&E^{(n)}_{W}(v_{1}\otimes \cdots\otimes v_{n}; w).
\end{eqnarray*}
This linear map has the $L(-1)$-derivative property,
the $L(0)$-conjugation property
and by Proposition \ref{correl-fn} 
is composable with $m$ vertex operators for any $m\in \Z_{+}$.
Moreover, 
let $f$ be a  homogeneous 
rational functions of degree $0$ in $z_{1}, \dots, z_{n}$ with the only possible 
poles at $z_{i}=z_{j}$, $i\ne j$, then 
$fE^{(n)}_{W;\;w}: V^{\otimes n}\to 
\widetilde{W}_{z_{1}, \dots, z_{n}}$ defined by 
\begin{eqnarray*}
\lefteqn{((fE^{(n)}_{W;\;w})(v_{1}\otimes \cdots \otimes v_{n}))(z_{1}, \dots, z_{n})}\nn
&&=f(z_{1}, \dots, z_{n})
(E^{(n)}_{W}(v_{1}\otimes \cdots\otimes v_{n}; w))(z_{1}, \dots, z_{n})
\end{eqnarray*}
for $v_{1}, \dots, v_{n}\in V$ has 
the $L(0)$-conjugation property
and is composable with $m$ vertex operators for any $m\in \Z_{+}$. In particular,
$fE^{(n)}_{V;\;\one}$ has 
the $L(0)$-conjugation property
and is composable with $m$ vertex operators for any $m\in \Z_{+}$.}
\end{expl}

Let $\Phi: V^{\otimes n}\to 
\widetilde{W}_{z_{1}, \dots, z_{n}}$
be composable with $m$ vertex operators. Then for  
$l_{1}, \ldots, l_{n} \in \Z_+$ such that $l_{1}+\cdots +l_{n}=m+n$,
$v_{1}, \dots, v_{m+n}\in V$ and $w\in W$, we have an element 
$$E(\Phi(E^{(l_{1})}_{V}(v_{1}\otimes \cdots \otimes v_{l_{1}}; \one)\otimes \cdots\otimes
E^{(l_{n})}_{V}(v_{l_{1}+\cdots +l_{n-1}+1}\otimes \cdots 
\otimes v_{l_{1}+\cdots +l_{n-1}+l_{n}}; \one)))$$
of 
$\widetilde{W}_{z_{1},  \dots, z_{m+n-1}}$ given by 
\begin{eqnarray*}
\lefteqn{\langle w', (E(\Phi(E^{(l_{1})}_{V}(v_{1}\otimes \cdots 
\otimes v_{l_{1}}; \one)\otimes}\nn 
&&\quad\quad  \cdots
\otimes E^{(l_{n})}_{V}(v_{l_{1}+\cdots +l_{n-1}+1}\otimes \cdots \otimes v_{l_{1}
+\cdots +l_{n-1}+l_{n}}; \one))))(z_{1}, \dots, z_{m+n})\rangle\nn
&&=R(\langle w', \Phi(E^{(l_{1})}_{V}(v_{1}\otimes \cdots 
\otimes v_{l_{1}}; \one)\otimes\nn
&&\quad\quad \cdots\otimes
E^{(l_{n})}_{V}(v_{l_{1}+\cdots +l_{n-1}+1}\otimes \cdots 
\otimes v_{l_{1}+\cdots +l_{n-1}+l_{n}}; \one))
(z_{1}, \dots, 
z_{m+n})\rangle)
\end{eqnarray*}
For $v_{1}, \dots, v_{m+n}\in V$, we have an element
$$E(E^{(m)}_{W}(v_{1}\otimes \cdots\otimes v_{m};
\Phi(v_{m+1}\otimes \cdots\otimes v_{m+n})))$$
of $\widetilde{W}_{z_{1}, \dots, 
z_{m+n}}$
given by 
\begin{eqnarray*}
\lefteqn{\langle w', (E(E^{(m)}_{W}(v_{1}\otimes \cdots\otimes v_{m};
\Phi(v_{m+1}\otimes \cdots\otimes v_{m+n})))(z_{1}, \dots, 
z_{m+n})\rangle}\nn
&&=R(\langle w', 
(E^{(m)}_{W}(v_{1}\otimes \cdots\otimes v_{m}; \nn
&&\quad\quad\quad\quad
(\Phi(v_{m+1}\otimes \cdots\otimes v_{m+n}))(z_{m+1}, \dots, z_{m+n}))
(z_{1}, \dots, z_{m})\rangle).
\end{eqnarray*}
Also for $v_{1}, \dots, v_{n+m}\in V$, since by Proposition \ref{wv<=>vw}, 
\begin{eqnarray}\label{wv-psi=vw-phi}
\lefteqn{\sum_{q\in \C}\langle w', 
(E^{W; (m)}_{WV}(P_{q}((\Phi(v_{1}\otimes \cdots\otimes v_{n}))(z_{1}, \dots, z_{n}))
; }\nn
&&\quad\quad\quad\quad\quad\quad\quad\quad
v_{n+1}\otimes \cdots\otimes v_{n+m}))
(z_{n+1}, \dots, z_{n+m})\rangle\nn
&&=\sum_{q\in \C}\langle w', 
(E^{(m)}_{W}(v_{n+1}\otimes \cdots\otimes v_{n+m}; \nn
&&\quad\quad\quad\quad
P_{q}((\Phi(v_{1}\otimes \cdots\otimes v_{n}))(z_{1}, \dots, z_{n}))))
(z_{n+1}, \dots, z_{n+m})\rangle,
\end{eqnarray}
the left-hand side of (\ref{wv-psi=vw-phi}) is absolutely convergent 
and can be analytically expended to a rational function 
in $z_{1}, \dots, z_{n}$ with the only possible 
poles at $z_{i}=z_{j}$, $i\ne j$
if and only if the same conclusions hold for 
the right-hand side of (\ref{wv-psi=vw-phi}). 
Since $\Phi$ is composable with $m$ vertex operators, 
the right-hand side is indeed absolutely convergent 
and can be analytically expended to a rational function
$z_{1}, \dots, z_{n}$ with the only possible 
poles at $z_{i}=z_{j}$, $i\ne j$. Thus the same conclusions
hold for the left-hand side. 
Denote the corresponding rational function by 
$$R(\langle w', 
(E^{W; (m)}_{WV}(\Phi(v_{1}\otimes \cdots\otimes v_{n})
; v_{n+1}\otimes \cdots\otimes v_{n+m}))
(z_{1}, \dots, z_{n+m})\rangle).$$
We obtain an element 
$$E(E^{W; (m)}_{WV}(\Phi(v_{1}\otimes \cdots\otimes v_{n})
; v_{n+1}\otimes \cdots\otimes v_{n+m}))$$
of $\widetilde{W}_{z_{1}, \dots, 
z_{n+m}}$
given by 
\begin{eqnarray*}
\lefteqn{\langle w', E(E^{W; (m)}_{WV}(\Phi(v_{1}\otimes \cdots\otimes v_{n})
; v_{n+1}\otimes \cdots\otimes v_{n+m}))(z_{1}, \dots, 
z_{n+m})\rangle}\nn
&&\!\!\!=R(\langle w', 
(E^{W; (m)}_{WV}(\Phi(v_{1}\otimes \cdots\otimes v_{n})
; v_{n+1}\otimes \cdots\otimes v_{n+m}))
(z_{1}, \dots, z_{n+m})\rangle).
\end{eqnarray*}
By Proposition \ref{wv<=>vw}, we have 
\begin{eqnarray}\label{e-wv<=>vw}
\lefteqn{E(E^{W; (m)}_{WV}(\Phi(v_{1}\otimes \cdots\otimes v_{n})
; v_{n+1}\otimes \cdots\otimes v_{n+m}))(z_{1}, \dots, 
z_{n+m})}\nn
&&\!\!\!=(E(E^{(m)}_{W}(v_{n+1}\otimes \cdots\otimes v_{n+m};
\Phi(v_{1}\otimes \cdots\otimes v_{n})))(z_{n+1}, \dots, z_{n+m}, z_{1},\cdots
z_{n}).\nn
\end{eqnarray}

We now define
$$\Phi\circ (E^{(l_{1})}_{V;\;\one}\otimes \cdots \otimes E^{(l_{n})}_{V;\;\one}): 
V^{\otimes m+n}\to 
\widetilde{W}_{z_{1},  \dots, z_{m+n}},$$
$$E^{(m)}_{W}\circ_{m+1}\Phi: V^{\otimes m+n}\to 
\widetilde{W}_{z_{1}, \dots, 
z_{m+n-1}}$$
and 
$$E^{W; (m)}_{WV}\circ_{0}\Phi: V^{\otimes m+n}\to 
\widetilde{W}_{z_{1}, \dots, 
z_{m+n-1}}$$
by
\begin{eqnarray*}
\lefteqn{(\Phi\circ (E^{(l_{1})}_{V;\;\one}\otimes \cdots \otimes 
E^{(l_{n})}_{V;\;\one}))(v_{1}\otimes \cdots \otimes v_{m+n-1})}\nn
&&=E(\Phi(E^{(l_{1})}_{V; \one}(v_{1}\otimes \cdots \otimes v_{l_{1}})\otimes \cdots
\nn 
&&\quad\quad\quad\quad\quad \otimes E^{(l_{n})}_{V; \one}
(v_{l_{1}+\cdots +l_{n-1}+1}\otimes \cdots 
\otimes v_{l_{1}+\cdots +l_{n-1}+l_{n}}))),
\end{eqnarray*}
\begin{eqnarray*}
\lefteqn{(E^{(m)}_{W}\circ_{m+1}\Phi)(v_{1}\otimes \cdots \otimes v_{m+n})}\nn
&&=E(E^{(m)}_{W}(v_{1}\otimes \cdots\otimes v_{m};
\Phi(v_{m+1}\otimes \cdots\otimes v_{m+n})))
\end{eqnarray*}
and 
\begin{eqnarray*}
\lefteqn{(E^{W; (m)}_{WV}\circ_{0}\Phi)(v_{1}\otimes \cdots \otimes v_{m+n})}\nn
&&=E(E^{W; (m)}_{WV}(\Phi(v_{1}\otimes \cdots\otimes v_{n})
; v_{n+1}\otimes \cdots\otimes v_{n+m})),
\end{eqnarray*}
respectively.
In the case that $l_{1}=\cdots=l_{i-1}=l_{i+1}=1$ and $l_{i}=m-n-1$ for some $i$,
for simplicity, we shall also use $\Phi\circ_{i} E^{(l_{i})}_{V;\;\one}$ to 
denote $\Phi\circ (E^{(l_{1})}_{V;\;\one}\otimes \cdots 
\otimes E^{(l_{n})}_{V;\;\one})$.

We define an action of $S_{n}$ on the space $\hom(V^{\otimes n}, 
\widetilde{W}_{z_{1}, \dots, z_{n}})$ of linear maps from 
$V^{\otimes n}$ to $\widetilde{W}_{z_{1}, \dots, z_{n}}$ by 
$$(\sigma(\Phi))(v_{1}\otimes \cdots\otimes v_{n})
=\sigma(\Phi(v_{\sigma(1)}\otimes \cdots\otimes v_{\sigma(n)}))$$
for $\sigma\in S_{n}$ and $v_{1}, \dots, v_{n}\in V$.

We shall use the notation $\sigma_{i_{1}, \dots, i_{n}}\in S_{n}$ to denote the 
the permutation given by 
$$\sigma_{i_{1}, \dots, i_{n}}(j)=i_{j}$$
for $j=1, \dots, n$.
We have
\begin{prop}
For $m\in \Z_{+}$, 
\begin{equation}\label{e-wv=simga-e-vw}
E^{W; (m)}_{WV}\circ_{0}\Phi
=\sigma_{n+1, \dots, n+m, 1,\dots,, n}(E^{(m)}_{W}\circ_{m+1}\Phi).
\end{equation}
\end{prop}
\pf
The equality (\ref{e-wv=simga-e-vw}) follows from
(\ref{e-wv<=>vw}) and the definition of the action of $S_{m+n}$
on $\widetilde{W}_{z_{1}, \dots, z_{m+n}}$.
\epfv

We also have:

\begin{prop}
The subspace of $\hom(V^{\otimes n}, 
\widetilde{W}_{z_{1}, \dots, z_{n}})$ consisting of linear maps
having
the $L(-1)$-derivative property, having the $L(0)$-conjugation property
or being composable with $m$ vertex operators is invariant under the 
action of $S_{n}$.
\end{prop}
\pf
This result follows directly from the definitions. 
\epfv

We know that compositions of maps are associative. But for maps 
whose compositions are defined using sums of 
absolutely convergent series as we have discussed above, 
even if all the compositions involved exist,
we might still not have associativity in general because iterated sums 
in different orders might not be equal to each other in general. However, when 
such compositions are analytic in some sense, associativity do hold. 
In particular, for the maps 
considered in this paper, we do have the following proposition that gives 
in particular some associativity results:

\begin{prop}\label{comp-assoc}
Let $\Phi: V^{\otimes n}\to 
\widetilde{W}_{z_{1}, \dots, z_{n}}$
be composable with $m$ vertex operators. Then we have:

\begin{enumerate}

\item For $p\le m$, $\Phi$ is 
composable with $p$ vertex operators and for 
$p, q\in \Z_{+}$ such that $p+q\le m$ and 
$l_{1}, \dots, l_{n} \in \Z_+$ such that $l_{1}+\cdots +l_{n}=p+n$,
$\Phi\circ (E^{(l_{1})}_{V;\;\one}\otimes 
\cdots \otimes E^{(l_{n})}_{V;\;\one})$ and $E^{(p)}_{W}\circ_{p+1}\Phi$
are
composable with $q$ vertex operators. 

\item For $p, q\in \Z_{+}$ such that $p+q\le m$,
$l_{1}, \dots, l_{n} \in \Z_+$ such that $l_{1}+\cdots +l_{n}=p+n$ and
$k_{1}, \dots, k_{p+n} \in \Z_+$ such that $k_{1}+\cdots +k_{p+n}=q+p+n$,
we have
\begin{eqnarray*}
&(\Phi\circ (E^{(l_{1})}_{V;\;\one}\otimes 
\cdots \otimes E^{(l_{n})}_{V;\;\one}))\circ 
(E^{(k_{1})}_{V;\;\one}\otimes 
\cdots \otimes E^{(k_{p+n})}_{V;\;\one})&\nn
&=\Phi\circ (E^{(k_{1}+\cdots +k_{l_{1}})}_{V;\;\one}\otimes 
\cdots \otimes E^{(k_{l_{1}+\cdots +l_{n-1}+1}+\cdots +k_{p+n})}_{V;\;\one}).&
\end{eqnarray*}

\item For $p, q\in \Z_{+}$ such that $p+q\le m$ and
$l_{1}, \dots, l_{n} \in \Z_+$ such that $l_{1}+\cdots +l_{n}=p+n$,
we have
$$E^{(q)}_{W}\circ_{q+1} (\Phi\circ (E^{(l_{1})}_{V;\;\one}\otimes 
\cdots \otimes E^{(l_{n})}_{V;\;\one}))
=(E^{(q)}_{W}\circ_{q+1} \Phi)\circ (E^{(l_{1})}_{V;\;\one}\otimes 
\cdots \otimes E^{(l_{n})}_{V;\;\one}).$$

\item For $p, q\in \Z_{+}$ such that $p+q\le m$, we have
$$E^{(p)}_{W}\circ_{p+1} (E^{(q)}_{W}\circ_{q+1}\Phi)
=E^{(p+q)}_{W}\circ_{p+q+1}\Phi.$$

\end{enumerate}
\end{prop}
\pf 
Conclusion 1 is clear from the definition. 

Let $v^{(i)}_{j}\in V$ for $i=1, \dots, p+n$, $j=1, \dots, 
k_{i}$. 
Then 
\begin{eqnarray}\label{comp-assoc-1}
\lefteqn{((\Phi\circ (E^{(l_{1})}_{V;\;\one}\otimes 
\cdots \otimes E^{(l_{n})}_{V;\;\one}))\circ(E^{(k_{1})}_{V;\;\one}\otimes 
\cdots \otimes E^{(k_{p+n})}_{V;\;\one}))(v^{(1)}_{1}\otimes 
\cdots \otimes v^{(p+n)}_{k_{p+n}})}\nn
&&=E((\Phi\circ (E^{(l_{1})}_{V;\;\one}\otimes 
\cdots \otimes E^{(l_{n})}_{V;\;\one}))(E^{(k_{1})}_{V;\;\one}(v^{1)}_{1}\otimes \cdots 
\otimes v^{(1)}_{k_{1}}) \otimes\nn
&&\quad\quad\quad\quad\quad \quad\quad \quad\quad \quad\quad \quad\quad 
\cdots \otimes E^{(k_{p+n})}_{V;\;\one}(v^{(p+n)}_{1}\otimes 
\cdots\otimes v^{(p+n)}_{k_{p+n}})))\nn
&&=E(\Phi(E^{(l_{1})}_{V;\;\one}(E^{(k_{1})}_{V;\;\one}(v^{(1)}_{1}\otimes \cdots 
\otimes v^{(1)}_{k_{1}}) \otimes
\cdots \otimes E^{(k_{l_{1}})}_{V;\;\one}(v^{(l_{1})}_{1}
\otimes \cdots 
\otimes v^{(l_{1})}_{k_{l_{1}}}))\otimes \nn
&&\quad\quad\quad
\cdots \otimes E^{(l_{n})}_{V;\;\one}(E^{(k_{l_{1}+\cdots +l_{n-1}+1})}_{V;\;\one}
(v^{(l_{1}+\cdots +l_{n-1}+1)}_{1}\otimes \cdots 
\otimes v^{(l_{1}+\cdots +l_{n-1}+1)}_{k_{l_{1}+\cdots +l_{n-1}+1}}) \otimes\nn
&&\quad\quad\quad\quad\quad \quad\quad
\cdots \otimes E^{(k_{p+n})}_{V;\;\one}(v^{(p+n)}_{1}\otimes 
\cdots\otimes v^{(p+n)}_{k_{p+n}})))).\nn
\end{eqnarray}
By Proposition \ref{correl-fn}, the right-hand side of (\ref{comp-assoc-1})
is equal to 
\begin{eqnarray}\label{comp-assoc-2}
\lefteqn{E(\Phi(E^{(k_{1}+\cdots +k_{l_{1}})}_{V;\;\one}
(v^{(1)}_{1}
\otimes \cdots 
\otimes  v^{(l_{1})}_{k_{l_{1}}})\otimes} \nn
&&\quad\quad\quad
\cdots \otimes E^{(k_{l_{1}+\cdots +l_{n-1}+1}+\cdots +k_{p+n})}_{V;\;\one}
(v^{(l_{1}+\cdots +l_{n-1}+1)}_{1}\otimes \cdots 
\otimes  v^{(p+n)}_{k_{p+n}})))\nn
&&=(\Phi\circ (E^{(k_{1}+\cdots +k_{l_{1}})}_{V;\;\one}\otimes 
\cdots \otimes E^{(k_{l_{1}+\cdots +l_{n-1}+1}+\cdots +k_{p+n})}_{V;\;\one}))
(v^{(1)}_{1}\otimes 
\cdots \otimes v^{(p+n)}_{k_{p+n}}).\nn
\end{eqnarray}
From (\ref{comp-assoc-1}) and (\ref{comp-assoc-2}), we obtain Conclusion 2.

Conclusions  3 and 4
can be proved similarly.
\epfv

\renewcommand{\theequation}{\thesection.\arabic{equation}}
\renewcommand{\thethm}{\thesection.\arabic{thm}}
\setcounter{equation}{0}
\setcounter{thm}{0}

\section{Chain complexes and cohomologies}

Let $V$ be a vertex operator algebra and $W$ a $V$-module. 
For $n\in \Z_{+}$, let $\widehat{C}_{0}^{n}(V, W)$ be the vector space of all 
linear maps from $V^{\otimes n}$ to $\widetilde{W}_{z_{1}, \dots, z_{n}}$ 
satisfying the $L(-1)$-derivative property and the $L(0)$-conjugation property. 
For $m, n\in \Z_{+}$, 
let $\widehat{C}_{m}^{n}(V, W)$ be the vector spaces of all 
linear maps from $V^{\otimes n}$ to $\widetilde{W}_{z_{1}, \dots, z_{n}}$
composable with $m$ vertex operators and satisfying the $L(-1)$-derivative
property and the $L(0)$-conjugation property. 
Also, let $\widehat{C}_{m}^{0}(V, W)=W$. Then we have 
$$\widehat{C}_{m}^{n}(V, W)\subset \widehat{C}_{m-1}^{n}(V, W)$$
for $m\in \Z_{+}$. Let 
$$\widehat{C}_{\infty}^{n}(V, W)=\bigcap_{m\in \N}\widehat{C}_{m}^{n}(V, W).$$
By Example \ref{correl-fns}, $\widehat{C}_{\infty}^{n}(V, V)$ is nonempty. 

For $n\in \N$ and $m\in \Z_{+}$, we define a 
coboundary operator 
$$\hat{\delta}^{n}_{m}: \widehat{C}_{m}^{n}(V, W)
\to \widehat{C}_{m-1}^{n+1}(V, W)$$
by 
\begin{eqnarray}\label{hat-delta}
\hat{\delta}^{n}_{m}(\Phi)&=&E^{(1)}_{W}\circ_{2} \Phi
+\sum_{i=1}^{n}(-1)^{i}\Phi\circ_{i} E^{(2)}_{V; \one}
 +(-1)^{n+1}
E^{W; (1)}_{WV}\circ_{0}
\Phi\nn
&=&E^{(1)}_{W}\circ_{2} \Phi
+\sum_{i=1}^{n}(-1)^{i}\Phi\circ_{i} E^{(2)}_{V; \one}
 +(-1)^{n+1}
\sigma_{n+1, 1, \dots, n}(E^{(1)}_{W}\circ_{2}
\Phi)\nn
\end{eqnarray}
for $\Phi \in \widehat{C}_{m}^{n}(V, W)$, where the second equality is 
obtained by using (\ref{e-wv=simga-e-vw}). Explicitly, 
for $v_{1}, \dots, v_{n+1}\in V$, $w'\in W'$ 
and $(z_{1}, \dots, z_{n+1})\in 
F_{n+1}\C$, 
\begin{eqnarray*}
\lefteqn{\langle w', ((\hat{\delta}^{n}_{m}(\Phi))(v_{1}\otimes \cdots\otimes v_{n+1}))
(z_{1}, \dots, z_{n+1})\rangle}\nn
&&=R(\langle w', Y_{W}(v_{1}, z_{1})(\Phi(v_{2}\otimes \cdots\otimes v_{n+1}))
(z_{2}, \dots, z_{n+1})\rangle)\nn
&&\quad +\sum_{i=1}^{n}(-1)^{i}R(\langle w', 
(\Phi(v_{1}\otimes \cdots \otimes v_{i-1} \nn
&&\quad\quad\quad\quad\quad\quad\quad
\otimes 
(Y_{V}(v_{i}, z_{i}-\zeta_{i})Y_{V}(v_{i+1}, z_{i+1}-\zeta_{i})\one)
\nn
&&\quad\quad\quad\quad\quad\quad\quad
\otimes v_{i+2}\otimes \cdots \otimes v_{n+1}))
(z_{1}, \dots, z_{i-1}, \zeta_{i}, z_{i+2}, \dots, z_{n+1})\rangle)\nn
&&\quad + (-1)^{n+1}R(\langle w', Y_{W}(v_{n+1}, z_{n+1})
(\Phi(v_{1}\otimes \cdots \otimes v_{n}))(z_{1}, \dots, z_{n})\rangle),
\end{eqnarray*}
which is in fact independent of $\zeta_{i}$. In particular, 
when we take $\zeta_{i}=z_{i+1}$ for $i=1, \dots, n$, we obtain
\begin{eqnarray*}
\lefteqn{\langle w', ((\hat{\delta}^{n}_{m}(\Phi))(v_{1}\otimes \cdots\otimes v_{n+1}))
(z_{1}, \dots, z_{n+1})\rangle}\nn
&&=R(\langle w', Y_{W}(v_{1}, z_{1})(\Phi(v_{2}\otimes \cdots\otimes v_{n+1}))
(z_{2}, \dots, z_{n+1})\rangle)\nn
&&\quad +\sum_{i=1}^{n}(-1)^{i}R(\langle w', 
(\Phi(v_{1}\otimes \cdots \otimes v_{i-1} \otimes 
Y_{V}(v_{i}, z_{i}-z_{i+1})v_{i+1}\nn
&&\quad\quad\quad\quad\quad\quad\quad\quad\quad 
\quad\quad\otimes \cdots \otimes v_{n+1}))
(z_{1}, \dots, z_{i-1}, z_{i+1}, \dots, z_{n+1})\rangle)\nn
&&\quad + (-1)^{n+1}R(\langle w', Y_{W}(v_{n+1}, z_{n+1})
(\Phi(v_{1}\otimes \cdots \otimes v_{n}))(z_{1}, \dots, z_{n})\rangle).
\end{eqnarray*}

By Proposition \ref{comp-assoc}, $\hat{\delta}^{n}_{m}(\Phi)$ is 
composable with $m-1$ vertex operators and has the $L(-1)$-derivative
property and the $L(0)$-conjugation property. So $\hat{\delta}^{n}_{m}(\Phi)\in 
\widehat{C}_{m-1}^{n+1}(V, W)$ and $\hat{\delta}^{n}_{m}$ is indeed a map whose image is in 
$\widehat{C}_{m-1}^{n+1}(V, W)$.

In the definition of  $\hat{\delta}^{n}_{m}$ above,
we require $m\in \Z_{+}$ so that each term in the right-hand sides of the
first and second equalities of (\ref{hat-delta}) is well defined. 
However, in the case $n=2$, there is a subspace 
of $\widehat{C}_{0}^{2}(V, W)$ 
containing $\widehat{C}_{m}^{2}(V, W)$ for all $m\in \Z_{+}$ such that 
$\hat{\delta}^{2}_{m}$ is still defined on this subspace. 

Let $\widehat{C}_{\frac{1}{2}}^{2}(V, W)$ be the subspace of $\widehat{C}_{0}^{2}(V, W)$ 
consisting of elements $\Phi$ such that for $v_{1}, v_{2}, v_{3}\in V$, $w'\in W'$, 
\begin{eqnarray*}
\lefteqn{\sum_{r\in \C}\big(\langle w', E^{(1)}_{W}(v_{1};
P_{r}((\Phi(v_{2}\otimes v_{3}))(z_{2}-\zeta, z_{3}-\zeta)))(z_{1}, \zeta)\rangle}\nn
&&\quad+\langle w', (\Phi(v_{1}\otimes P_{r}((E^{(2)}_{V}(v_{2}\otimes v_{3}; \one))
(z_{2}-\zeta, z_{3}-\zeta))))
(z_{1}, \zeta)\rangle\big)
\end{eqnarray*}
and 
\begin{eqnarray*}
\lefteqn{\sum_{r\in \C}\big(\langle w', 
(\Phi(P_{r}((E^{(2)}_{V}(v_{1}\otimes v_{2}; \one))(z_{1}-\zeta, z_{2}-\zeta))
\otimes v_{3}))
(\zeta, z_{3})\rangle}\nn
&&\quad +\langle w', 
E^{W; (1)}_{WV}(P_{r}((\Phi(v_{1}\otimes v_{2}))(z_{1}-\zeta, z_{2}-\zeta));
v_{3}))(\zeta, z_{3})\rangle\big)
\end{eqnarray*}
are absolutely convergent in the regions $|z_{1}-\zeta|>|z_{2}-\zeta|, |z_{2}-\zeta|>0$ and 
$|\zeta-z_{3}|>|z_{1}-\zeta|, |z_{2}-\zeta|>0$, respectively, 
and can be analytically extended to 
rational functions in $z_{1}$ and $z_{2}$ with the only possible poles at
$z_{1}, z_{2}=0$ and $z_{1}=z_{2}$. Note that here we do not require 
the individual series 
\begin{eqnarray*}
&{\displaystyle \sum_{r\in \C}\big(\langle w', E^{(1)}_{W}(v_{1};
P_{r}((\Phi(v_{2}\otimes v_{3}))(z_{2}-\zeta, z_{3}-\zeta)))(z_{1}, \zeta)\rangle,}&\\
&{\displaystyle \sum_{r\in \C}\langle w', 
(\Phi(v_{1}\otimes P_{r}((E^{(2)}_{V}(v_{2}\otimes v_{3}; \one))
(z_{2}-\zeta, z_{3}-\zeta))))
(z_{1}, \zeta)\rangle,}&\\
&{\displaystyle \sum_{r\in \C}\langle w', 
(\Phi(P_{r}((E^{(2)}_{V}(v_{1}\otimes v_{2}; \one))(z_{1}-\zeta, z_{2}-\zeta))
\otimes v_{3}))
(\zeta, z_{3})\rangle,}&\\
&{\displaystyle \sum_{r\in \C}\langle w', 
E^{W; (1)}_{WV}(P_{r}((\Phi(v_{1}\otimes v_{2}))(z_{1}-\zeta, z_{2}-\zeta));
v_{3}))(\zeta, z_{3})\rangle}&
\end{eqnarray*}
to be absolutely convergent. 
We denote  the corresponding rational functions 
by 
\begin{eqnarray*}
\lefteqn{R(\langle w', (E^{(1)}_{W}(v_{1};
\Phi(v_{2}\otimes v_{3}))(z_{1}, z_{2}, z_{3})\rangle}\nn
&&\quad +\langle w', (\Phi(v_{1}\otimes E^{(2)}_{V}(v_{2}\otimes v_{3}; \one)))
(z_{1}, z_{2}, z_{3})\rangle)
\end{eqnarray*}
and 
\begin{eqnarray*}
\lefteqn{R(\langle w', 
(\Phi(E^{(2)}_{V}(v_{1}\otimes v_{2}; \one))
\otimes v_{3}))(z_{1}, z_{2}, z_{3})\rangle}\nn
&&\quad +\langle w', 
(E^{W; (1)}_{WV}(\Phi(v_{1}\otimes v_{2}); v_{3}))
(z_{1}, z_{2}, z_{3})\rangle).
\end{eqnarray*}
Clearly, 
$\widehat{C}_{m}^{2}(V, W)\subset \widehat{C}_{\frac{1}{2}}^{2}(V, W)$
for $m\in \Z_{+}$. We define a coboundary operator 
$$\hat{\delta}^{2}_{\frac{1}{2}}: \widehat{C}_{\frac{1}{2}}^{2}(V, W)
\to \widehat{C}_{0}^{3}(V, W)$$
by 
\begin{eqnarray*}
\lefteqn{\langle w', ((\hat{\delta}^{2}_{\frac{1}{2}}(\Phi))
(v_{1}\otimes v_{2} \otimes v_{3}))(z_{1}, z_{2}, z_{3})\rangle}\nn
&&=R(\langle w', (E^{(1)}_{W}(v_{1};
\Phi(v_{2}\otimes v_{3}))(z_{1}, z_{2}, z_{3})\rangle\nn
&&\quad \quad +\langle w', (\Phi(v_{1}\otimes E^{(2)}_{V}(v_{2}\otimes v_{3}; \one)))
(z_{1}, z_{2}, z_{3})\rangle)\nn
&&\quad -R(\langle w', 
(\Phi(E^{(2)}_{V}(v_{1}\otimes v_{2}; \one))
\otimes v_{3}))(z_{1}, z_{2}, z_{3})\rangle\nn
&&\quad \quad +\langle w', 
(E^{W; (1)}_{WV}(\Phi(v_{1}\otimes v_{2}); v_{3}))
(z_{1}, z_{2}, z_{3})\rangle)
\end{eqnarray*}
for $w'\in W'$, 
$\Phi\in \widehat{C}_{\frac{1}{2}}^{2}(V, W)$,
$v_{1}, v_{2}, v_{3}\in V$ and $(z_{1}, z_{2}, z_{3})\in F_{3}\C$. 

%\begin{eqnarray*}
%\lefteqn{(\hat{\delta}_{n}(\Phi))(v_{1}\otimes \cdots\otimes v_{n+1})}\nn
%&&=E(E^{(1)}_{W}(v_{1}; \Phi(v_{2}\otimes \cdots\otimes v_{n+1})))\nn
%&&\quad +\sum_{i=1}^{n}(-1)^{i}E(\Phi(v_{1}\otimes \cdots\otimes 
%v_{i-1}\nn
%&& \quad\quad\quad\quad\quad\quad\quad\quad
% \otimes E^{(2)}_{V}(v_{i}\otimes v_{i+1})
%\otimes v_{i+2}\otimes \cdots\otimes 
%v_{n+1}))\nn
%&&\quad +(-1)^{n+1}
%\sigma_{n+1, 1, \dots, n}(E(E^{(1)}_{W}(v_{n+1}; \Phi(v_{1}\otimes \cdots\otimes v_{n})))),
%\end{eqnarray*}

\begin{prop}\label{delta-square}
For $n\in \N$ and $m\in \Z_{+}+1$, $\hat{\delta}^{n+1}_{m-1}\circ \hat{\delta}^{n}_{m}=0$.
We also have $\hat{\delta}^{2}_{\frac{1}{2}}\circ \hat{\delta}^{1}_{2}=0$.
\end{prop}
\pf
Let $\Phi\in \widehat{C}^{n}(V, W)$. 
Then 
\begin{eqnarray}\label{delta-square-1}
\lefteqn{(\hat{\delta}^{n+1}_{m-1}\circ \hat{\delta}^{n}_{m})(\Phi)
}\nn
&&=E^{(1)}_{W}\circ_{2}(\hat{\delta}^{n}_{m}(\Phi))
 +\sum_{i=1}^{n+1}(-1)^{i}(\hat{\delta}^{n}_{m}(\Phi))\circ_{i}
\otimes E^{(2)}_{V; \one}\nn
&&\quad 
+(-1)^{n+2}
\sigma_{n+2, 1, \dots, n+1}(E^{(1)}_{W}\circ_{2}
(\hat{\delta}^{n}_{m}(\Phi)))\nn
&&=E^{(1)}_{W}\circ_{2}(E^{(1)}_{W}\circ_{2}\Phi)
+\sum_{j=1}^{n}(-1)^{j}E^{(1)}_{W}\circ_{2}(\Phi\circ_{j+1}
E^{(2)}_{V; \one})\nn
&&\quad +(-1)^{n+1}E^{(1)}_{W}\circ_{2}
(\sigma_{n+1, 1, \dots, n}(E_{W}^{(1)}\circ_{2}\Phi))
\nn
&&\quad -(E^{(1)}_{W}\circ_{2}
\Phi)\circ_{1} E^{(2)}_{V; \one}
+\sum_{i=2}^{n+1}(-1)^{i}(E^{(1)}_{W}\circ_{2}\Phi)\circ_{i}
E^{(2)}_{V; \one}\nn
&&\quad +\sum_{i=1}^{n+1}(-1)^{i}\sum_{j=1}^{i-2}(-1)^{j}
(\Phi\circ_{j}E^{(2)}_{V; \one})\circ_{i} 
E^{(2)}_{V; \one}\nn
&&\quad 
+\sum_{i=2}^{n+1}(-1)^{i}(-1)^{i-1}(\Phi\circ_{i-1}
E^{(2)}_{V; \one})\circ_{i}E^{(2)}_{V; \one}\nn
&&\quad +\sum_{i=1}^{n}(-1)^{i}(-1)^{i}(\Phi\circ_{i}E^{(2)}_{V; \one})
\circ_{i}E^{(2)}_{V; \one}\nn
&&\quad +\sum_{i=1}^{n+1}(-1)^{i}\sum_{j=i+2}^{n+1}(-1)^{j-1}
(\Phi\circ_{j}E^{(2)}_{V; \one})\circ_{i}E^{(2)}_{V; \one}\nn
&&\quad +(-1)^{n+1}\sum_{i=1}^{n}(-1)^{i}(\sigma_{n+1, 1, \dots, n}
(E^{(1)}_{W}\circ_{2}\Phi))\circ_{i}E^{(2)}_{V; \one}\nn
&&\quad + (\sigma_{n+1, 1, \dots, n}
(E^{(1)}_{W}\circ_{2}\Phi))\circ_{n+1}E^{(2)}_{V; \one}\nn
&&\quad +(-1)^{n+2}
\sigma_{n+2, 1, \dots, n+1}(E^{(1)}_{W}\circ_{2}(E^{(1)}_{W}\circ_{2}
\Phi))
\nn
&&\quad +(-1)^{n+2}\sum_{j=1}^{n}(-1)^{j}
\sigma_{n+2, 1, \dots, n+1}(E^{(1)}_{W}\circ_{2}(\Phi\circ_{j}E^{(2)}_{V; \one}))\nn
&&\quad -\sigma_{n+2, 1, \dots, n+1}(E^{(1)}_{W}\circ_{2}
\sigma_{n+1, 1, \dots, n}
(E^{(1)}_{W}\circ_{2}\Phi)).
\end{eqnarray}
We now prove that in the right-hand side of (\ref{delta-square-1}),
(i) the first  and the fourth terms, (ii) the second and the 
fifth terms, (iii) the third and the twelfth terms, (iv) the sixth and the ninth term,
(v) the seventh and the eighth terms, (vi) the tenth and the thirteenth terms, 
(vii) the eleventh and fourteenth terms cancel with each other, and thus the 
right-hand side of (\ref{delta-square}) is equal to $0$, proving the
proposition. In the proofs of these cancellations below, we actually also have to 
switch the order of absolutely convergent iterated sums. But just as in the proof of 
Proposition \ref{comp-assoc}, since all the iterated sums are absolutely convergent 
to rational functions, the corresponding multisums are all absolutely convergent
and thus all iterated sums are equal. Because of this general fact, we shall 
omit the discussion of the orders of the iterated sums.

(i) The first  and the fourth terms: 
For $w'\in W'$, $v_{1}, \dots, v_{n+2}\in V$ and $(z_{1}, \dots, z_{n+2})\in F_{n+2}\C$, 
applying the first term to $v_{1}\otimes \cdots
\otimes v_{n+2}$, evaluating the result at $(z_{1}, \dots, z_{n+2})$ and then 
pairing the result with $w'$, we obtain
\begin{eqnarray}\label{d-2=0-1}
\lefteqn{\langle w', (E(E^{(1)}_{W}(v_{1}; 
E(E^{(1)}_{W}(v_{2}; \Phi(v_{3}\otimes \cdots\otimes v_{n+2}))))))(z_{1}, \dots, 
z_{n+1})\rangle}\nn
&&=R(\langle w', (E^{(1)}_{W}(v_{1}; 
E(E^{(1)}_{W}(v_{2}; \Phi(v_{3}\otimes \cdots\otimes v_{n+2})))))(z_{1}, \dots, 
z_{n+1})\rangle)\nn
&&=R(\langle w', Y_{W}(v_{1}, z_{1})
E(E^{(1)}_{W}(v_{2}; \Phi(v_{3}\otimes \cdots\otimes v_{n+2})))(z_{1}, \dots, 
z_{n+1})\rangle)\nn
&&=R(\langle w', Y_{W}(v_{1}, z_{1})Y_{W}(v_{2}, z_{2})
(\Phi(v_{3}\otimes \cdots\otimes v_{n+2}))(z_{3}, \dots, z_{n+2})
\rangle)\nn
&&=R(\langle w', Y_{W}(Y_{V}(v_{1}, z_{1}-z_{2})v_{2}, z_{2})
(\Phi(v_{3}\otimes \cdots\otimes v_{n+2}))(z_{3}, \dots, z_{n+2})
\rangle),\nn
\end{eqnarray}
where in the last step, we have used the associativity of the module $W$. 
On the other hand, applying the fourth term to $v_{1}\otimes \cdots
\otimes v_{n+2}$, evaluating the result at $(z_{1}, \dots, z_{n+2})$ and then 
pairing the result with $w'$, we obtain
\begin{eqnarray}\label{d-2=0-2}
\lefteqn{\langle w', -((E^{(1)}_{W}\circ_{2}\Phi)(E^{(2)}_{V}(v_{1}\otimes v_{2}; \one)
\otimes v_{3}\otimes \cdots\otimes v_{n+2}))(z_{1}, \dots, z_{n+2})\rangle}\nn
&&=-\langle w', (E(E^{(1)}_{W}(E^{(2)}_{V}(v_{1}\otimes v_{2}; \one); 
\Phi(v_{3}\otimes \cdots\otimes v_{n+2})))))(z_{1}, \dots, z_{n+2})\rangle\nn
&&=-R(\langle w', (E^{(1)}_{W}(E^{(2)}_{V}(v_{1}\otimes v_{2}; \one); 
\Phi(v_{3}\otimes \cdots\otimes v_{n+2})))(z_{1}, \dots, z_{n+2})\rangle)\nn
&&=-R(\langle w', Y_{W}(Y_{V}(v_{1}, z_{1}-z_{2})v_{2}, z_{2}) 
(\Phi(v_{3}\otimes \cdots\otimes v_{n+2}))(z_{3}, \dots, z_{n+2})\rangle).\nn
\end{eqnarray}
Since the the right-hand sides of (\ref{d-2=0-1}) and (\ref{d-2=0-2}) 
differ only by a sign for $w'\in W'$, $v_{1}, \dots, v_{n+2}\in V$
and $(z_{1}, \dots, z_{n+2})\in F_{n+2}\C$, these two terms indeed cancel with each other. 

(ii) The second and the 
fifth terms: By Proposition \ref{comp-assoc}, these two terms
differ only by a sign and thus cancel with each other. 

(iii) The third and the twelfth terms:
For $w'\in W'$, $v_{1}, \dots, v_{n+2}\in V$ and $(z_{1}, \dots, z_{n+2})\in F_{n+2}\C$,  
applying the third term to $v_{1}\otimes \cdots
\otimes v_{n+2}$, evaluating the result at $(z_{1}, \dots, z_{n+2})$ and then 
pairing the result with $w'$, we obtain
\begin{eqnarray}\label{d-2=0-3}
\lefteqn{\langle w', (-1)^{n+1}(E(E^{(1)}_{W}(v_{1}; }\nn
&& \quad\quad\quad
\sigma_{n+1, 1, \dots, n}(E(E_{W}(v_{n+2}; \Phi(v_{2}\otimes \cdots\otimes v_{n+1})))))))
(z_{1}, \dots, z_{n+2})\rangle\nn
&&=(-1)^{n+1}R(\langle w', (E^{(1)}_{W}(v_{1}; \nn
&& \quad\quad\quad
\sigma_{n+1, 1, \dots, n}(E(E_{W}(v_{n+2}; \Phi(v_{2}\otimes \cdots\otimes v_{n+1}))))))
(z_{1}, \dots, z_{n+2})\rangle)\nn
&&=(-1)^{n+1}R(\langle w', Y_{W}(v_{1}, z_{1})\cdot \nn
&& \quad\quad\quad\cdot
\sigma_{n+1, 1, \dots, n}(E(E_{W}(v_{n+2}; \Phi(v_{2}\otimes \cdots\otimes v_{n+1}))))
(z_{2}, \dots, z_{n+2})\rangle)\nn
&&=(-1)^{n+1}R(\langle w', Y_{W}(v_{1}, z_{1}) \cdot\nn
&& \quad\quad\quad\cdot
(E(E_{W}(v_{n+2}; \Phi(v_{2}\otimes \cdots\otimes v_{n+1}))))
(z_{n+2}, z_{2}\dots, z_{n+1})\rangle)\nn
&&=(-1)^{n+1}R(\langle w', Y_{W}(v_{1}, z_{1})Y_{W}(v_{n+2}, z_{n+2}) \cdot\nn
&& \quad\quad\quad\quad\quad\quad\quad\quad\quad\quad\quad
\cdot(\Phi(v_{2}\otimes \cdots\otimes v_{n+1}))
(z_{2}\dots, z_{n+1})\rangle)\nn
&&=(-1)^{n+1}R(\langle w', Y_{W}(v_{n+2}, z_{n+2}) Y_{W}(v_{1}, z_{1})\cdot\nn
&& \quad\quad\quad\quad\quad\quad\quad\quad\quad\quad\quad
\cdot(\Phi(v_{2}\otimes \cdots\otimes v_{n+1}))
(z_{2}\dots, z_{n+1})\rangle),
\end{eqnarray}
where in the last step, we have used the commutativity of the module $W$. 
On the other hand, applying the twelfth term to $v_{1}\otimes \cdots
\otimes v_{n+2}$, evaluating the result at $(z_{1}, \dots, z_{n+2})$ and then 
pairing the result with $w'$, we obtain
\begin{eqnarray}\label{d-2=0-4}
\lefteqn{\langle w', (-1)^{n+2}
\sigma_{n+2, 1, \dots, n+1}(E(E^{(1)}_{W}(v_{n+2}; E(E^{(1)}_{W}(v_{1}; }\nn
&& \quad\quad\quad\quad\quad\quad\quad\quad\quad\quad\quad\quad
\Phi(v_{2}\otimes \cdots\otimes v_{n+1}))))))(z_{1}, \dots, z_{n+2})\rangle\nn
&&=(-1)^{n+2}R(\langle w', (E^{(1)}_{W}(v_{n+2}; E(E^{(1)}_{W}(v_{1}; \nn
&& \quad\quad\quad\quad\quad\quad\quad\quad\quad\quad\quad\quad
\Phi(v_{2}\otimes \cdots\otimes v_{n+1})))))(z_{n+2}, z_{1}, \dots, z_{n+1})\rangle)\nn
&&=(-1)^{n+2}R(\langle w', Y_{W}(v_{n+2}, z_{n+2})Y_{W}(v_{1}, z_{1})\cdot\nn
&& \quad\quad\quad\quad\quad\quad\quad\quad\quad\quad\quad\quad\cdot
(\Phi(v_{2}\otimes \cdots\otimes v_{n+1}))(z_{2}, \dots, z_{n+1})\rangle).
\end{eqnarray}
Since the right-hand sides of (\ref{d-2=0-3}) and (\ref{d-2=0-4}) 
differ only by a sign for $w'\in W'$, $v_{1}, \dots, v_{n+2}\in V$
and $(z_{1}, \dots, z_{n+2})\in F_{n+2}\C$, these two terms indeed cancel with each other. 

(iv) The sixth and the ninth terms: 
By Proposition \ref{comp-assoc}, when $i\ne j$,
we have 
$$(\Phi\circ_{j}E^{(2)}_{V; \one})\circ_{i}E^{(2)}_{V; \one}=\Phi\circ 
(E_{V; \one}^{(l_{1})}\otimes \cdots\otimes E_{V; \one}^{(l_{n})}),$$
where $l_{k}=1$ if $k\ne i, j$ and $l_{k}=2$ if $k=i$ or $k=j$. 
Since for any $a_{ij}\in C^{n+2}(V, W)$,
$$\sum_{j=1}^{n+1}\sum_{i=1}^{j-2}a_{ij}=\sum_{i=1}^{n+1}\sum_{j=i+2}^{n+1}a_{ij},$$
we see that these two terms differ by a sign and 
thus cancel with each other.

(v) The seventh and the eighth terms:
By Proposition \ref{comp-assoc}, we have 
$$(\Phi\circ_{i-1}
E^{(2)}_{V; \one})\circ_{i}E^{(2)}_{V; \one}=\Phi\circ_{i-1}E^{(3)}_{V; \one}$$
for $i=2, \dots, n+1$ and 
and
$$(\Phi\circ_{i}E^{(2)}_{V; \one})
\circ_{i}E^{(2)}_{V; \one}=\Phi\circ_{i}E^{(3)}_{V; \one}$$
for $i=1, \dots, n$.
Thus these two terms differ by a sign and 
cancel with each other.

(vi) The tenth and the thirteenth terms: Since $i$ or $j$ are less than or equal to 
$n$, by definition, we have 
$$(\sigma_{n+1, 1, \dots, n}(E_{W}^{(1)}\circ_{2}\Phi))\circ_{i}E_{V; \one}^{(2)}
=\sigma_{n+2, 1, \dots, n+1}(E_{W}^{(1)}\circ_{2}(\Phi\circ_{i}E_{V; \one}^{(2)}).$$
Then by Proposition \ref{comp-assoc}, these two terms differ by a sign and 
thus cancel with each other. 

(vii) The eleventh and the fourteenth terms: 
For $w'\in W'$, $v_{1}, \dots, v_{n+2}\in V$ and $(z_{1}, \dots, z_{n+2})\in F_{n+2}\C$,  
applying the eleventh term to $v_{1}\otimes \cdots
\otimes v_{n+2}$, evaluating the result at $(z_{1}, \dots, z_{n+2})$ and then 
pairing the result with $w'$, we obtain
\begin{eqnarray}\label{d-2=0-9}
\lefteqn{\langle w', 
((\sigma_{n+1, 1, \dots, n}(E^{(1)}_{W}\circ_{2}\Phi))(v_{1}\otimes
\cdots\otimes v_{n}}\nn
&&\quad\quad\quad\quad\quad\quad\quad\quad\quad\quad\quad\quad
\otimes E_{V}^{(2)}(v_{n+1}\otimes v_{n+2}; \one)))(z_{1}, 
\dots, z_{n+2})\rangle\nn
&&=\langle w', ((E^{(1)}_{W}\circ_{2}\Phi)(E^{(2)}_{V}(v_{n+1}\otimes v_{n+2}; \one) \nn
&& \quad\quad\quad\quad\quad\quad\quad\quad\quad\quad\quad\quad
\otimes v_{1}\otimes \cdots\otimes v_{n}))(z_{n+1}, z_{n+2}, z_{1}, 
\dots, z_{n})\rangle\nn
&&=\langle w', 
(E(E^{(1)}_{W}(E^{(2)}_{V}(v_{n+1}\otimes v_{n+2}; \one); \nn
&& \quad\quad\quad\quad\quad\quad\quad\quad\quad\quad\quad\quad
\Phi(v_{1}\otimes \cdots\otimes v_{n}))))(z_{n+1}, z_{n+2}, z_{1}, 
\dots, z_{n})\rangle\nn
&&=R(\langle w', 
(E^{(1)}_{W}(E^{(2)}_{V}(v_{n+1}\otimes v_{n+2}; \one); \nn
&& \quad\quad\quad\quad\quad\quad\quad\quad\quad\quad\quad\quad
 \Phi(v_{1}\otimes \cdots\otimes v_{n})))(z_{n+1}, z_{n+2}, z_{1}, 
\dots, z_{n})\rangle)\nn
&&=R(\langle w', Y_{W}(Y_{V}(v_{n+1}, z_{n+1}-z_{n+2})v_{n+2}, z_{n+2})\cdot\nn
&& \quad\quad\quad\quad\quad\quad\quad\quad\quad\quad\quad\quad\cdot
(\Phi(v_{1}\otimes \cdots\otimes v_{n}))(z_{1}, 
\dots, z_{n})\rangle)\nn
&&=R(\langle w', Y_{W}(v_{n+1}, z_{n+1})Y_{W}(v_{n+2}, z_{n+2})
(\Phi(v_{1}\otimes \cdots\otimes v_{n}))(z_{1}, 
\dots, z_{n})\rangle)\nn
&&=R(\langle w', Y_{W}(v_{n+2}, z_{n+2})Y_{W}(v_{n+1}, z_{n+1})
(\Phi(v_{1}\otimes \cdots\otimes v_{n}))(z_{1}, 
\dots, z_{n})\rangle),\nn
\end{eqnarray}
where in the last two steps, we have used the associativity and 
commutativity of the $V$-module $W$. 
On the other hand, applying the fourteenth term to $v_{1}\otimes \cdots
\otimes v_{n+2}$, evaluating the result at $(z_{1}, \dots, z_{n+2})$ and then 
pairing the result with $w'$, we obtain
\begin{eqnarray}\label{d-2=0-10}
\lefteqn{\langle w',((E^{(1)}_{W}\circ_{2}
\sigma_{n+1, 1, \dots, n}
(E^{(1)}_{W}\circ_{2}\Phi))(v_{n+2}\otimes v_{1}\otimes }\nn
&& \quad\quad\quad\quad\quad\quad\quad\quad\quad\quad\quad\quad\quad\quad\quad
 \cdots \otimes 
v_{n+1}))(z_{n+2}, z_{1}, \dots, z_{n+1})\rangle \nn
&&=-\langle w', (E(E^{(1)}_{W}(v_{n+2}; \nn
&& \quad\quad\quad
((\sigma_{n+1, 1, \dots, n}(E^{(1)}_{W}\circ_{2}\Phi))(v_{1}\otimes 
\cdots\otimes v_{n+1})))))
(z_{n+2}, z_{1}, \dots, z_{n+1})\rangle\nn
&&=-R(\langle w', Y_{W}(v_{n+2}, z_{n+2})\cdot\nn
&& \quad\quad\quad\cdot
((\sigma_{n+1, 1, \dots, n}(E^{(1)}_{W}\circ_{2}\Phi))(v_{1}\otimes 
\cdots\otimes v_{n+1}))
(z_{1}, \dots, z_{n+1})\rangle)\nn
&&=-R(\langle w', Y_{W}(v_{n+2}, z_{n+2})\cdot\nn
&& \quad\quad\quad\cdot
((E^{(1)}_{W}\circ_{2}\Phi)(v_{n+1}\otimes v_{1}
\cdots\otimes v_{n}))
(z_{n+1}, z_{1} \dots, z_{n})\rangle)\nn
&&=-R(\langle w', Y_{W}(v_{n+2}, z_{n+2})\cdot\nn
&& \quad\quad\quad\cdot
(E(E^{(1)}_{W}(v_{n+1}; 
\Phi(v_{1}\otimes \cdots\otimes v_{n}))))
(z_{n+1}, z_{1}, \dots, z_{n})\rangle)\nn
&&=-R(\langle w', Y_{W}(v_{n+2}, z_{n+2})Y_{W}(v_{n+1}, z_{n+1})
(\Phi(v_{1}\otimes \cdots\otimes v_{n}))
(z_{1}, \dots, z_{n})\rangle).\nn
\end{eqnarray}
Since the right-hand sides of (\ref{d-2=0-9}) and (\ref{d-2=0-10}) 
differ only by a sign for $w'\in W'$, $v_{1}, \dots, v_{n+2}\in V$
and $(z_{1}, \dots, z_{n+2})\in F_{n+2}\C$, these two terms indeed cancel with each other.

The second conclusion follows from the first conclusion. In fact since 
$$\hat{\delta}_{2}^{1}(\widehat{C}_{2}^{1}(V, W))\subset 
\widehat{C}_{1}^{2}(V, W)\subset 
\widehat{C}_{\frac{1}{2}}^{2}(V, W),$$
$$\hat{\delta}^{2}_{\frac{1}{2}}\circ \hat{\delta}^{1}_{2}
=\hat{\delta}^{2}_{1}\circ \hat{\delta}^{1}_{2}
=0.$$
\epfv

By Proposition \ref{delta-square},  we have complexes
\begin{equation}\label{hat-complex}
0\longrightarrow \widehat{C}_{m}^{0}(V, W)
\stackrel{\hat{\delta}^{0}_{m}}{\longrightarrow}
\widehat{C}_{m-1}^{1}(V, W)
\stackrel{\hat{\delta}^{1}_{m-1}}{\longrightarrow}\cdots
\stackrel{\hat{\delta}^{m-1}_{1}}{\longrightarrow}
\widehat{C}_{0}^{m}(V, W)\longrightarrow 0
\end{equation}
for $m\in \N$ and 
\begin{equation}\label{hat-complex-half}
0\longrightarrow \widehat{C}^{0}_{3}(V, W)
\stackrel{\hat{\delta}_{3}^{0}}{\longrightarrow}
\widehat{C}^{1}_{2}(V, W)
\stackrel{\hat{\delta}_{2}^{1}}{\longrightarrow}\widehat{C}_{\frac{1}{2}}^{2}(V, W)
\stackrel{\hat{\delta}_{\frac{1}{2}}^{2}}{\longrightarrow}
\widehat{C}_{0}^{3}(V, W)\longrightarrow 0,
\end{equation}
where the first and last arrows are the trivial 
embeddings and projections.
But these complexes are not yet the chain complexes for $V$. We have to consider 
certain subcomplexes of these complexes. 
To define these subcomplexes, we need to use
shuffles. 

For $n\in \N$ and $1\le p \le n-1$, let $J_{n; p}$ be the set of elements of 
$S_{n}$ which preserve the order of the first $p$ numbers and the order of the last 
$n-p$ numbers, that is,
$$J_{n, p}=\{\sigma\in S_{n}\;|\;\sigma(1)<\cdots <\sigma(p),\;
\sigma(p+1)<\cdots <\sigma(n)\}.$$
Elements of $J_{n; p}$ are called {\it shuffles}. Let $J_{n; p}^{-1}=\{\sigma\;|\;
\sigma\in J_{n; p}\}$.
For $m, n\in \N$ and or $m=\frac{1}{2}$, $n=2$, let $C_{m}^{n}(V, W)$ 
be the subspace of 
$\widehat{C}_{m}^{n}(V, W)$ consisting of maps $\Phi$ such that
$$\sum_{\sigma\in J_{n; p}^{-1}}(-1)^{|\sigma|}
\sigma(\Phi(v_{\sigma(1)}\otimes \cdots \otimes v_{\sigma(n)}))=0.$$
We also let 
$$C_{\infty}^{n}(V, W)=\bigcap_{m\in \N}C_{m}^{n}(V, W).$$

\begin{thm}\label{delta-for-c}
For $n\in \N$ and $m\in \Z_{+}$, 
$\hat{\delta}^{n}_{m}(C_{m}^{n}(V, W))\subset C_{m-1}^{n+1}(V, W)$.
Also $\hat{\delta}^{2}_{\frac{1}{2}}(C_{\frac{1}{2}}^{2}(V, W))\subset C_{0}^{3}(V, W)$.
\end{thm}
\pf
Let $\Phi\in C_{m}^{n}(V, W)$. We need to prove 
$\hat{\delta}^{n}_{m}(\Phi)\in C_{m-1}^{n+1}(V, W)$.
By definition, we have 
\begin{eqnarray}\label{delta-for-c-1}
\lefteqn{\sum_{\sigma\in J_{n+1; p}^{-1}}(-1)^{|\sigma|}
\sigma(\hat{\delta}^{n}_{m}(\Phi))}\nn
&&=\sum_{\sigma\in J_{n+1; p}^{-1}}(-1)^{|\sigma|}
\sigma(E^{(1)}_{W}\circ_{2}\Phi)\nn
&&\quad +\sum_{\sigma\in J_{n+1; p}^{-1}}(-1)^{|\sigma|}
\sum_{i=1}^{n}(-1)^{i}\sigma(\Phi\circ_{i}E^{(2)}_{V; \one})\nn
&&\quad +(-1)^{n+1}\sum_{\sigma\in J_{n+1; p}^{-1}}(-1)^{|\sigma|}
 \sigma(\sigma_{n+1, 1 \dots, n}(E^{(1)}_{W}\circ_{2}\Phi)).
\end{eqnarray}

Note that for any  $\sigma\in J_{n+1; p}^{-1}$, 
$\sigma(1)$ is either $1$ or $p+1$. 
So the first term in the right-hand side of (\ref{delta-for-c-1}) is equal to 
\begin{equation}\label{delta-for-c-2}
\sum_{\mbox{\scriptsize $\begin{array}{c}
\sigma\in J_{n+1; p}^{-1}\\ \sigma(1)=1\end{array}$}}(-1)^{|\sigma|}
\sigma(E^{(1)}_{W}\circ_{2}\Phi)
+\sum_{\mbox{\scriptsize $\begin{array}{c}\sigma\in J_{n+1; p}^{-1}\\
\sigma(1)=p+1\end{array}$}}(-1)^{|\sigma|}
\sigma(E^{(1)}_{W}\circ_{2}\Phi).
\end{equation}
The subset $\{\sigma\in J_{n+1; p}^{-1}\;|\;\sigma(1)=1\}$ 
of $S_{n+1}$ is the image of $J_{n; p}^{-1}$ under the embedding from $S_{n}$ to $S_{n+1}$ 
given by mapping an element of $S_{n}$ to an element of $S_{n+1}$ permuting 
only the last $n$ numbers. Since $\Phi\in C_{m}^{n}(V, W)$, the first term in 
(\ref{delta-for-c-2}) is equal to $0$. Similarly, the subset
$\{\sigma\in J_{n+1; p}^{-1}\;|\;\sigma(1)=p+1\}$ 
of $S_{n+1}$ is the image of $J_{n; p}^{-1}$ under the embedding from $S_{n}$ to $S_{n+1}$ 
given by mapping an element of $S_{n}$ to an element of $S_{n+1}$ permuting 
only the numbers $1, \dots, p, p+2, \dots, n+1$. So the second term in
(\ref{delta-for-c-2}) is also equal to $0$. Thus the first term in
the right-hand side of (\ref{delta-for-c-1}) is equal to $0$. Similarly, with 
$n+1$ and $p$ playing the role of $1$ and $p+1$ in the argument above,
the last term in the right-hand side of (\ref{delta-for-c-1}) is also 
equal to $0$.

The second term in the right-hand side of (\ref{delta-for-c-1}) is equal to 
\begin{eqnarray}\label{delta-for-c-3}
\lefteqn{\sum_{i=1}^{p}
(-1)^{i}\sum_{\mbox{\scriptsize $\begin{array}{c}
\sigma\in J_{n+1; p}^{-1}\\ 1\le \sigma(i),\; \sigma(i+1)\le p\end{array}$}}(-1)^{|\sigma|}
\sigma(\Phi\circ_{i}E^{(2)}_{V; \one})}\nn
&&+\sum_{i=1}^{p}
(-1)^{i}\sum_{\mbox{\scriptsize $\begin{array}{c}
\sigma\in J_{n+1; p}^{-1}\\ p+1\le \sigma(i),\; \sigma(i+1)\le n\end{array}$}}(-1)^{|\sigma|}
\sigma(\Phi\circ_{i}E^{(2)}_{V; \one})\nn
&&+\sum_{i=p+1}^{n}
(-1)^{i}\sum_{\mbox{\scriptsize $\begin{array}{c}
\sigma\in J_{n+1; p}^{-1}\\ 1\le \sigma(i)\le p\\ p+1\le\sigma(i+1)\le n\end{array}$}}
(-1)^{|\sigma|}
\sigma(\Phi\circ_{i}E^{(2)}_{V; \one})\nn
&&+\sum_{i=p+1}^{n}
(-1)^{i}\sum_{\mbox{\scriptsize $\begin{array}{c}
\sigma\in J_{n+1; p}^{-1}\\ p+1\le \sigma(i)\le n\\ 1\le
\sigma(i+1)\le p\end{array}$}}(-1)^{|\sigma|}
\sigma(\Phi\circ_{i}E^{(2)}_{V; \one}).
\end{eqnarray}
For $\sigma
\in J_{n+1; p}^{-1}$, if $1\le \sigma(i), \sigma(i+1)\le p$ or 
$p+1\le \sigma(i), \sigma(i+1)\le n$,
then $\sigma(i+1)=\sigma(i)+1$. In these cases, the permutation of 
the numbers $1, \dots, i-1, i+1, \dots, n$ induced from $\sigma^{-1}$ 
is in $J_{n; p-1}^{-1}$ or $J_{n; p}^{-1}$ and 
every element of $J_{n; p-1}^{-1}$ or $J_{n; p}^{-1}$
is obtained uniquely in this way. Since $\Phi\in C_{m}^{n}(V, W)$, the first two terms 
in (\ref{delta-for-c-3}) are equal to $0$. 
Let $\sigma \in J_{n+1; p}^{-1}$ satisfying $p+1\le \sigma(i)\le n,\; 1\le
\sigma(i+1)\le p$. Let $\sigma_{i+1, i}$ be the transposition exchanging $i$ and $i+1$.
Then $\tau=\sigma\circ\sigma_{i+1, i}$ is an element of $J_{n+1; p}^{-1}$
satisfying $1\le \tau(i)\le p,\; p+1\le\tau(i+1)\le n$. Moreover, 
$|\sigma\circ\sigma_{i+1, i}|=|\sigma|+1$. Also, by the commutativity of $V$,
$$(\sigma\circ\sigma_{i+1, i})(\Phi\circ_{i}E^{(2)}_{V; \one})
=\sigma(\Phi\circ_{i}E^{(2)}_{V; \one}).$$
Thus the third term and the fourth term
in (\ref{delta-for-c-3}) cancel with each other. 
The calculations above show that 
(\ref{delta-for-c-3}) is equal to $0$. 

Now we have proved that the right-hand side of 
(\ref{delta-for-c-1}) is equal to $0$. By (\ref{delta-for-c-1}),
$\hat{\delta}^{n}_{m}(\Phi)\in C_{m-1}^{n+1}(V, W)$. 

Let $\Phi\in C_{\frac{1}{2}}^{2}(V, W)$. Note that $J_{3; 2}=\sigma_{3, 1, 2}(J_{3; 1})$.
To prove $\hat{\delta}^{2}_{\frac{1}{2}}(\Phi)\in C_{0}^{3}(V, W)$,
we need only prove 
\begin{equation}\label{delta-for-c-3.5}
\sum_{\sigma\in J_{3; 1}^{-1}}(-1)^{|\sigma|}\sigma(\hat{\delta}^{2}_{\frac{1}{2}}(\Phi))=0.
\end{equation}
By definition, for $v_{1}, v_{2}, v_{3}\in V$, $w'\in W'$, 
$(z_{1}, z_{2}, z_{3})\in F_{3}\C$ and $\zeta\in \C$ such that 
$(z_{1}-\zeta, z_{2}-\zeta), (z_{2}-\zeta, z_{3}-\zeta),
(z_{1}-\zeta, z_{3}-\zeta), (z_{1}, \zeta), (z_{2}, \zeta),
(z_{3}, \zeta)\in F_{2}\C$, we have
\begin{eqnarray}\label{delta-for-c-4}
\lefteqn{\sum_{\sigma\in J_{3; 1}^{-1}}(-1)^{|\sigma|}
\langle w', \sigma(\hat{\delta}^{2}_{\frac{1}{2}}(\Phi))
(v_{1}\otimes v_{2} \otimes v_{3}))(z_{1}, z_{2}, z_{3})\rangle}\nn
&&=\sum_{\sigma\in J_{3; 1}^{-1}}(-1)^{|\sigma|}
\langle w', (\hat{\delta}^{2}_{\frac{1}{2}}(\Phi))
(v_{\sigma(1)}\otimes v_{\sigma(2)} \otimes v_{\sigma(3)}))(z_{\sigma(1)}, 
z_{\sigma(2)}, z_{\sigma(3)})\rangle\nn
&&=\langle w', (\hat{\delta}^{2}_{\frac{1}{2}}(\Phi))
(v_{1}\otimes v_{2} \otimes v_{3}))(z_{1}, 
z_{2}, z_{3})\rangle\nn
&&\quad -\langle w', (\hat{\delta}^{2}_{\frac{1}{2}}(\Phi))
(v_{2}\otimes v_{1} \otimes v_{3}))(z_{2}, 
z_{1}, z_{3})\rangle\nn
&&\quad +\langle w', (\hat{\delta}^{2}_{\frac{1}{2}}(\Phi))
(v_{2}\otimes v_{3} \otimes v_{1}))(z_{2}, 
z_{3}, z_{1})\rangle\nn
&&=R(\langle w', (E^{(1)}_{W}(v_{1};
\Phi(v_{2}\otimes v_{3})))(z_{1}, z_{2}, z_{3})\rangle)\nn
&&\quad \quad\quad+\langle w', (\Phi(v_{1}\otimes 
(E^{(2)}_{V}(v_{2}\otimes v_{3}; \one)))(z_{1}, z_{2}, z_{3})\rangle)\nn
&&\quad -R(\langle w', 
(\Phi(E^{(2)}_{V}(v_{1}\otimes v_{2}; \one)
\otimes v_{3}))(z_{1}, z_{2}, z_{3})\rangle\nn
&&\quad\quad\quad +\langle w', 
(E^{W; (1)}_{WV}(\Phi(v_{1}\otimes v_{2}); v_{3}))
(z_{1}, z_{2}, z_{3})\rangle)\nn
&&\quad -
R(\langle w', (E^{(1)}_{W}(v_{2};
\Phi(v_{1}\otimes v_{3})))(z_{1}, z_{2}, z_{3})\rangle)\nn
&&\quad \quad\quad+\langle w', (\Phi(v_{2}\otimes 
E^{(2)}_{V}(v_{1}\otimes v_{3}; \one))(z_{1}, z_{2}, z_{3})\rangle)\nn
&&\quad +R(\langle w', 
(\Phi(E^{(2)}_{V}(v_{2}\otimes v_{1}; \one)
\otimes v_{3}))(z_{1}, z_{2}, z_{3})\rangle\nn
&&\quad\quad\quad +\langle w', 
(E^{W; (1)}_{WV}(\Phi(v_{2}\otimes v_{1}); v_{3}))
(z_{1}, z_{2}, z_{3})\rangle)\nn
&&\quad +
R(\langle w', (E^{(1)}_{W}(v_{2};
\Phi(v_{3}\otimes v_{1})))(z_{1}, z_{2}, z_{3})\rangle)\nn
&&\quad \quad\quad+\langle w', (\Phi(v_{2}\otimes 
(E^{(2)}_{V}(v_{3}\otimes v_{1}; \one)))(z_{1}, z_{2}, z_{3})\rangle)\nn
&&\quad -R(\langle w', 
(\Phi(E^{(2)}_{V}(v_{2}\otimes v_{3}; \one)
\otimes v_{1}))(z_{1}, z_{2}, z_{3})\rangle\nn
&&\quad\quad\quad -\langle w', 
(E^{W; (1)}_{WV}(\Phi(v_{2}\otimes v_{3}); v_{1}))(z_{1}, z_{2}, z_{3})\rangle)
\end{eqnarray}

Since $\Phi\in C_{\frac{1}{2}}^{2}(V, W)$ and $E^{(2)}_{V; \one}\in C_{\frac{1}{2}}^{2}(V, V)$,
we have 
\begin{eqnarray*}
(\Phi(u_{1}\otimes u_{2}))(\zeta_{1}, \zeta_{2})
&=&(\Phi(u_{2}\otimes u_{1}))(\zeta_{2}, \zeta_{1}),\\
(E^{(2)}_{V}(u_{1}\otimes u_{2}; \one))(\zeta_{1}, \zeta_{2})
&=&(E^{(2)}_{V}(u_{2}\otimes u_{1}; \one))(\zeta_{2}, \zeta_{1})
\end{eqnarray*}
for $u_{1}, u_{2}\in V$ and $(\zeta_{1}, \zeta_{2})\in F_{2}\C$. 
Also 
$$E^{(1)}_{W}(u; w)=E^{W; (1)}_{WV}(w; u)$$
for $u\in V$ and $w\in W$. From these formulas, we see that the first and 
sixth terms, the second and fourth terms and the third and fifth terms in 
(\ref{delta-for-c-4}) cancel with each other, proving (\ref{delta-for-c-3.5}).
\epfv

Let 
$\delta_{m}^{n}=\hat{\delta}_{m}^{n}|_{C_{m}^{n}(V, W)}$ for $m\in \Z_{+}$, 
$n\in \N$ or 
$m=\frac{1}{2}$, $n=2$. 
By Theorem \ref{delta-for-c},
the image of $\delta^{n}_{m}$ is in $C_{m-1}^{n+1}(V, W)$. So we obtain a 
linear map
$$\delta^{n}_{m}: C_{m}^{n}(V, W)\to C_{m-1}^{n+1}(V, W)$$
for each pair $m\in \Z_{+}$, $n\in \N$ or 
$$\delta^{2}_{\frac{1}{2}}: C_{\frac{1}{2}}^{2}(V, W)\to C_{0}^{3}(V, W).$$ 
For $m\in \Z_{+}$, we have a subcomplex
\begin{equation}\label{complex}
0\longrightarrow C_{m}^{0}(V, W)
\stackrel{\delta^{0}_{m}}{\longrightarrow}
C_{m-1}^{1}(V, W)
\stackrel{\delta^{1}_{m-1}}{\longrightarrow}\cdots
\stackrel{\delta^{m-1}_{1}}{\longrightarrow}
C_{0}^{m}(V, W)\longrightarrow 0
\end{equation}
of the complex (\ref{hat-complex}) and 
also a subcomplex 
\begin{equation}\label{complex-half}
0\longrightarrow C_{3}^{0}(V, W)
\stackrel{\delta_{3}^{0}}{\longrightarrow}
C_{2}^{1}(V, W)
\stackrel{\delta_{2}^{1}}{\longrightarrow}C_{\frac{1}{2}}^{2}(V, W)
\stackrel{\delta_{\frac{1}{2}}^{2}}{\longrightarrow}
C_{0}^{3}(V, W)\longrightarrow 0
\end{equation}
of the complex (\ref{hat-complex-half}).

Since $C_{\infty}^{n}(V, W)\subset C_{m}^{n}(V, W)$
for any $m\in \Z_{+}$ and
$C_{m_{2}}^{n}(V, W)\subset C_{m_{1}}^{n}(V, W)$,
for $m_{1}, m_{2}\in \Z_{+}$ satisfying $m_{1}\le m_{2}$, 
$\delta_{m}^{n}\lbar_{C_{\infty}^{n}(V, W)}$ is independent of 
$m$. Let 
$$\delta_{\infty}^{n}=\delta_{m}^{n}\lbar_{C_{\infty}^{n}(V, W)}: 
C_{\infty}^{n}(V, W)
\to C_{\infty}^{n+1}(V, W).$$
We obtain  a complex 
\begin{equation}\label{complex-infty}
0\longrightarrow C_{\infty}^{0}(V, W)
\stackrel{\delta_{\infty}^{0}}{\longrightarrow}
C_{\infty}^{1}(V, W)
\stackrel{\delta_{\infty}^{1}}{\longrightarrow}C_{\infty}^{2}(V, W)
\stackrel{\delta_{\infty}^{2}}{\longrightarrow}
\cdots.
\end{equation}

Using the complexes (\ref{complex}), (\ref{complex-half})
and (\ref{complex-infty}), 
we now introduce the cohomology spaces
of $V$.

\begin{defn}
{\rm For $m\in \Z_{+}$ and $n\in \N$, we define the 
{\it $n$-th cohomology $H^{n}_{m}(V, W)$ of $V$ with coefficient in 
$W$ and composable with 
$m$ vertex operators}
to be 
$$H_{m}^{n}(V, W)=\ker \delta^{n}_{m}/\mbox{\rm im}\; \delta^{n-1}_{m+1}.$$
We also define 
$$H^{2}_{\frac{1}{2}}(V, W)
=\ker \delta^{2}_{\frac{1}{2}}/\mbox{\rm im}\; \delta_{2}^{1}$$
and 
$$H_{\infty}^{n}(V, W)
=\ker \delta^{n}_{\infty}/\mbox{\rm im}\; \delta^{n-1}_{\infty}$$
for $n\in \N$.}
\end{defn}

For $n\in \N$ and $m_{1}, m_{2}\in \Z_{+}$, since 
$C_{m_{2}}^{n}(V, W)\subset C_{m_{1}}^{n}(V, W)$,
we have 
\begin{eqnarray*}
\delta_{m_{2}+1}^{n-1}(C_{m_{2}+1}^{n-1}(V, W))&\subset&
 \delta_{m_{1}+1}^{n-1}(C_{m_{2}+1}^{n-1}(V, W))\cap C_{m_{2}}^{n}(V, W),\\
\ker \delta_{m_{2}}^{n}&\subset& \ker \delta_{m_{1}}^{n}
\cap C_{m_{2}}^{n}(V, W).
\end{eqnarray*}
Thus we have an injective linear map
$f_{m_{1}m_{2}}: H_{m_{2}}^{n}(V, W)\to H_{m_{1}}^{n}(V, W)$
given by 
$$f_{m_{1}m_{2}}(\Phi+\ker \delta_{m_{2}}^{n})
=\Phi+\ker \delta_{m_{1}}^{n}.$$

\begin{prop}
For $n\in \N$, $(H_{m}^{n}(V, W), f_{m_{1}m_{2}})$
is an inverse system. Moreover, their inverse limits
are linearly isomorphic to $H_{\infty}^{n}(V, W)$ for $n\in \N$. 
\end{prop}
\pf
This follows straightforwardly from the definitions. 
\epfv

\begin{prop}
Let $V$ be a grading-restricted vertex algebra and 
$W$ a generalized $V$-module. Then $H_{m}^{0}(V, W)=W$ 
for any $m\in \Z_{+}$. 
\end{prop}
\pf
By definition, for any $m\in \Z_{+}$, $C_{m}^{0}(V, W)=W$ and 
$\delta_{m}^{0} (C_{m}^{0}(V, W))=0$.
So $H_{m}^{0}(V, W)=C_{m}^{0}(V, W)=W$.
\epfv

We shall discuss the first and second cohomologies in \cite{Hu3}.

\noindent {\small \sc Department of Mathematics, Rutgers University,
110 Frelinghuysen Rd., Piscataway, NJ 08854-8019 (permanent address)}
\vspace{1em}

\noindent  {\it and}
\vspace{1em}

\noindent {\small \sc Beijing International Center for Mathematical
Research, Peking University, Beijing 100871, China}
\vspace{1em}

\noindent {\em E-mail address}: yzhuang@math.rutgers.edu

\end{document}